\documentclass[reqno,12pt]{amsart}
\usepackage{amssymb}
\usepackage{amsfonts}
\usepackage{amsmath}
\usepackage{graphicx}
\usepackage{amscd,amsthm}

\newtheorem{theorem}{Theorem}
\theoremstyle{plain}
\newtheorem{acknowledgement}{Acknowledgement}

\newtheorem{definition}{Definition}

\newtheorem{remark}{Remark}

\numberwithin{equation}{section}

\input{tcilatex}

\begin{document}
\author{}
\title{}
\maketitle

\begin{center}
\thispagestyle{empty} \pagestyle{myheadings} 
\markboth{\bf Yilmaz Simsek
}{\bf New Generating Functions of the Stirling numbers, Frobenius-Euler and Related polynomials}

\textbf{{\Large Generating functions for generalized \textbf{Stirling type
numbers, }Array type polynomials,\textbf{\ }Eulerian type polynomials and
their applications}}

\bigskip

\textbf{Yilmaz Simsek}\\[0pt]

Department of Mathematics, Faculty of Science University of Akdeniz TR-07058
Antalya, Turkey ysimsek@akdeniz.edu.tr\\[0pt]

\bigskip

\textbf{{\large {Abstract}}}\medskip
\end{center}

\begin{quotation}
The first aim of this paper is to construct new generating functions for the
generalized $\lambda $-Stirling type numbers of the second kind, generalized
array type polynomials and generalized Eulerian type polynomials and
numbers, attached to Dirichlet character. We derive various functional
equations and differential equations using these generating functions. The
second aim is provide a novel approach to deriving identities including
multiplication formulas and recurrence relations\textit{\ }for these numbers
and polynomials using these functional equations and differential equations.
Furthermore, by applying $p$-adic Volkenborn integral and Laplace transform,
we derive some new identities for the generalized $\lambda $-Stirling type
numbers of the second kind, the generalized array type polynomials and the
generalized Eulerian type polynomials. We also give many applications
related to the class of these polynomials and numbers.
\end{quotation}

\bigskip

\noindent \textbf{2010 Mathematics Subject Classification.} 12D10, 11B68,
11S40, 11S80, 26C05, 26C10, 30B40, 30C15.

\bigskip

\noindent \textbf{Key Words.} Bernoulli polynomials; Euler polynomials;
Apostol Bernoulli polynomials; generalized Frobenius Euler polynomials;
Normalized Polynomials; Array polynomials; Stirling numbers of the second
kind; $p$-adic Volkenborn integral; generating function; functional
equation; Laplace transform.

\section{Introduction, Definitions and Preliminaries}

Throughout this paper, we use the following standard notations:

$\mathbb{N}=\{1,2,3,$\ldots $\}$, $\mathbb{N}_{0}=\{0,1,2,3,$\ldots $\}=%
\mathbb{N}\cup \{0\}$ and $\mathbb{Z}^{-}=\{-1,-2,-3,$\ldots $\}$. Here, $%
\mathbb{Z}$ denotes the set of integers, $\mathbb{R}$ denotes the set of
real numbers and $\mathbb{C}$ denotes the set of complex numbers. We assume
that $\ln (z)$ denotes the principal branch of the multi-valued function $%
\ln (z)$ with the imaginary part $\Im \left( \ln (z)\right) $ constrained by%
\begin{equation*}
-\pi <\Im \left( \ln (z)\right) \leq \pi .
\end{equation*}%
Furthermore, 
\begin{equation*}
0^{n}=\left\{ 
\begin{array}{cc}
1 & n=0 \\ 
&  \\ 
0 & n\in \mathbb{N},%
\end{array}%
\right. 
\end{equation*}%
\begin{equation*}
\left( 
\begin{array}{c}
x \\ 
v%
\end{array}%
\right) =\frac{x(x-1)\cdots (x-v+1)}{v!}
\end{equation*}%
and%
\begin{equation*}
\left\{ z\right\} _{0}=1\text{ and }\left\{ z\right\}
_{j}=\dprod\limits_{d=0}^{j-1}(z-d),
\end{equation*}%
where $j\in \mathbb{N}$ and $z\in \mathbb{C}$ cf. (\cite{Comtet}, \cite%
{LuoSrivatava2010}).

The generating functions have various applications in many branches of
Mathematics and Mathematical Physics. These functions are defined by linear
polynomials, differential relations, globally referred to as \textit{%
functional equations}. The functional equations arise in well-defined
combinatorial contexts and they lead systematically to well-defined classes
of functions (cf. see, for detail, \cite{Flajolet}). Although, in the
literature, one can find extensive investigations related to the generating
functions for the Bernoulli, Euler and Genocchi numbers and polynomials and
also their generalizations, the $\lambda $-Stirling numbers of the second
kind, the array polynomials and the Eulerian polynomials, related to
nonnegative real parameters, have not been studied yet. Therefore, Section
2, Section 3 and Section 4 of this paper deal with new classes of generating
functions which are related to generalized $\lambda $-Stirling type numbers
of the second kind, generalized array type polynomials and generalized
Eulerian polynomials, respectively. By using these generating functions, we
derive many functional equations and differential equations. By using these
equations, we investigate and introduce fundamental properties and many new
identities for the generalized $\lambda $-Stirling type numbers of the
second kind, the generalized array type polynomials and the generalized
Eulerian type polynomials and numbers. We also derive multiplication
formulas and recurrence relations for these numbers and polynomials.

The remainder of this study is organized as follows:

In section 5, we derive new identities related to the generalized Bernoulli
polynomials, the generalized Eulerian type polynomials, generalized $\lambda 
$-Stirling type numbers and the generalized array polynomials.

In section 6, we give relations between generalized Bernoulli polynomials
and generalized array polynomials.

In section 7, We give an application of the Laplace transform to the
generating functions for the generalized Bernoulli polynomials and the
generalized array type polynomials.

In section 8, by using the bosonic and the fermionic $p$-adic integral on $%
\mathbb{Z}_{p}$, we find some new identities related to the Bernoulli
polynomials, the generalized Eulerian type polynomials and Stirling numbers.

\section{Generating Function for generalized $\protect\lambda $-Stirling
type numbers of the second kind}

The Stirling numbers are used in combinatorics, in number theory, in
discrete probability distributions for finding higher order moments, etc.
The Stirling number of the second kind, denoted by $S(n,k)$, is the number
of ways to partition a set of $n$ objects into $k$ groups. These numbers
occur in combinatorics and in the theory of partitions.

In this section, we construct a new generating function, related to
nonnegative real parameters, for the generalized $\lambda $-Stirling type
numbers of the second kind. We derive some elementary properties including
recurrence relations of these numbers. The following definition provides a
natural generalization and unification of the $\lambda $-Stirling numbers of
the second kind:

\begin{definition}
Let $a$,$~b\in \mathbb{R}^{+}$ ($a\neq b$), $\lambda \in \mathbb{C}$ and $%
v\in \mathbb{N}_{0}$. The generalized $\lambda $-Stirling type numbers of
the second kind $\mathcal{S}(n,v;a,b;\lambda )$\ are defined by means of the
following generating function:%
\begin{equation}
f_{S,v}(t;a,b;\lambda )=\frac{\left( \lambda b^{t}-a^{t}\right) ^{v}}{v!}%
=\sum_{n=0}^{\infty }\mathcal{S}(n,v;a,b;\lambda )\frac{t^{n}}{n!}.
\label{s1}
\end{equation}
\end{definition}

\begin{remark}
By setting $a=1$ and $b=e$ in (\ref{s1}), we have the $\lambda $-Stirling
numbers of the second kind%
\begin{equation*}
\mathcal{S}(n,v;1,e;\lambda )=S(n,v;\lambda )
\end{equation*}%
which are defined by means of the following generating function:%
\begin{equation*}
\frac{\left( \lambda e^{t}-1\right) ^{v}}{v!}=\sum_{n=0}^{\infty
}S(n,v;\lambda )\frac{t^{n}}{n!},
\end{equation*}%
cf. (\cite{LuoSrivatava2010}, \cite{Srivastava2011}). Substituting $\lambda
=1$ into above equation, we have the Stirling numbers of the second kind%
\begin{equation*}
S(n,v;1)=S(n,v),
\end{equation*}%
cf. (\cite{Comtet}, \cite{LuoSrivatava2010}, \cite{Srivastava2011}). These
numbers have the following well known properties:%
\begin{equation*}
S(n,0)=\delta _{n,0},
\end{equation*}%
\begin{equation*}
S(n,1)=S(n,n)=1
\end{equation*}%
and%
\begin{equation*}
S(n,n-1)=\left( 
\begin{array}{c}
n \\ 
2%
\end{array}%
\right) ,
\end{equation*}%
where $\delta _{n,0}$ denotes the Kronecker symbol (see \cite{Comtet}, \cite%
{LuoSrivatava2010}, \cite{Srivastava2011}).
\end{remark}

By using (\ref{s1}), we obtain the following theorem:

\begin{theorem}
\label{Theorem STnumber}%
\begin{equation}
\mathcal{S}(n,v;a,b;\lambda )=\frac{1}{v!}\sum_{j=0}^{v}(-1)^{j}\left( 
\begin{array}{c}
v \\ 
j%
\end{array}%
\right) \lambda ^{v-j}\left( j\ln a+(v-j)\ln b\right) ^{n}  \label{as1}
\end{equation}%
and%
\begin{equation}
\mathcal{S}(n,v;a,b;\lambda )=\frac{1}{v!}\sum_{j=0}^{v}(-1)^{v-j}\left( 
\begin{array}{c}
v \\ 
j%
\end{array}%
\right) \lambda ^{j}\left( j\ln b+(v-j)\ln a\right) ^{n}.  \label{as1a}
\end{equation}
\end{theorem}

\begin{proof}
By using (\ref{s1}) and the binomial theorem, we can easily arrive at the
desired results.
\end{proof}

By using the formula (\ref{as1}), we can compute some values of the numbers $%
\mathcal{S}(n,v;a,b;\lambda )$ as follows:%
\begin{equation*}
\mathcal{S}(0,0;a,b;\lambda )=1,
\end{equation*}%
\begin{equation*}
\mathcal{S}(0,0;a,b;\lambda )=1,
\end{equation*}%
\begin{equation*}
\mathcal{S}(1,0;a,b;\lambda )=0,
\end{equation*}%
\begin{equation*}
\mathcal{S}(1,1;a,b;\lambda )=\ln \left( \frac{b^{\lambda }}{a}\right) ,
\end{equation*}%
\begin{equation*}
\mathcal{S}(2,0;a,b;\lambda )=0,
\end{equation*}%
\begin{equation*}
\mathcal{S}(2,1;a,b;\lambda )=\lambda \left( \ln b\right) ^{2}-\left( \ln
a\right) ^{2},
\end{equation*}%
\begin{equation*}
\mathcal{S}(2,2;a,b;\lambda )=\frac{\lambda ^{2}}{2}\left( \ln b^{2}\right)
^{2}-\lambda \ln \left( ab\right) +\left( \ln a^{2}\right) ^{2},
\end{equation*}

\begin{equation*}
\mathcal{S}(3,0;a,b;\lambda )=0,
\end{equation*}%
\begin{equation*}
\mathcal{S}(3,1;a,b;\lambda )=\lambda \left( \ln b\right) ^{3}-\left( \ln
a\right) ^{3},
\end{equation*}%
\begin{equation*}
\mathcal{S}(0,v;a,b;\lambda )=\frac{\left( \lambda -1\right) ^{v}}{v!},
\end{equation*}%
\begin{equation*}
\mathcal{S}(n,0;a,b;\lambda )=\delta _{n,0}
\end{equation*}%
and%
\begin{equation*}
\mathcal{S}(n,1;a,b;\lambda )=\lambda \left( \ln b\right) ^{n}-\left( \ln
a\right) ^{n}.
\end{equation*}

\begin{remark}
By setting $a=1$ and $b=e$ in the assertions (\ref{as1}) of Theorem \ref%
{Theorem STnumber}, we have the following result:%
\begin{equation*}
S(n,v;\lambda )=\frac{1}{v!}\sum_{j=0}^{v}\left( 
\begin{array}{c}
v \\ 
j%
\end{array}%
\right) \lambda ^{v-j}(-1)^{j}\left( v-j\right) ^{n}.
\end{equation*}%
The above relation has been studied by Srivastava \cite{Srivastava2011} and
Luo \cite{LuoSrivatava2010}. By setting $\lambda =1$ in the above equation,
we have the following result:%
\begin{equation*}
S(n,v;\lambda )=\frac{1}{v!}\sum_{j=0}^{v}\left( 
\begin{array}{c}
v \\ 
j%
\end{array}%
\right) (-1)^{j}\left( v-j\right) ^{n}
\end{equation*}%
cf. (\cite{AgohDilcher}, \cite{cagic}, \cite{Carlitz}, \cite{Carlitz1953G}, 
\cite{Comtet}, \cite{T. Kim}, \cite{LuoSrivatava2010}, \cite{SimsekSpringer}%
, \cite{YsimsekStirling}, \cite{Srivastava2011}, \cite{SrivastawaGargeSC}).
\end{remark}

By differentiating both sides of equation (\ref{s1}) with respect to the
variable $t$, we obtain the following\textit{\ }differential equations:%
\begin{equation*}
\frac{\partial }{\partial t}f_{S,v}(t;a,b;\lambda )=\left( \lambda (\ln
b)b^{t}-(\ln a)a^{t}\right) f_{S,v-1}(t;a,b;\lambda )
\end{equation*}%
or%
\begin{equation}
\frac{\partial }{\partial t}f_{S,v}(t;a,b;\lambda )=v\ln
(b)f_{S,v}(t;a,b;\lambda )+\ln \left( \frac{b}{a}\right)
a^{t}f_{S,v-1}(t;a,b;\lambda ).  \label{s1a}
\end{equation}

By using equations (\ref{s1}) and (\ref{s1a}), we obtain recurrence
relations for the generalized $\lambda $-Stirling type numbers of the second
kind by the following theorem:

\begin{theorem}
\label{TE2} Let $n,v\in \mathbb{N}$.%
\begin{equation}
\mathcal{S}(n,v;a,b;\lambda )=\sum_{j=0}^{n-1}\left( 
\begin{array}{c}
n-1 \\ 
j%
\end{array}%
\right) \mathcal{S}(j,v-1;a,b;\lambda )\left( \lambda \left( \ln (b)\right)
^{n-j}-\left( \ln (a)\right) ^{n-j}\right) .  \label{s4}
\end{equation}%
or%
\begin{eqnarray*}
\mathcal{S}(n,v;a,b;\lambda ) &=&v\ln (b)\mathcal{S}(n-1,v;a,b;\lambda ) \\
&&+\ln \left( \frac{b}{a}\right) \sum_{j=0}^{n-1}\left( 
\begin{array}{c}
n-1 \\ 
j%
\end{array}%
\right) \mathcal{S}(j,v-1;a,b;\lambda )\left( \ln (a)\right) ^{n-1-j}.
\end{eqnarray*}
\end{theorem}

\begin{remark}
By setting $a=1$ and $b=e$, Theorem \ref{TE2} yields the corresponding
results which are proven by Luo and Srivastava \cite[Theorem 11]%
{LuoSrivatava2010}. Substituting $a=\lambda =1$ and $b=e$ into Theorem \ref%
{TE2}, we obtain the following known results:%
\begin{equation*}
S(n,v)=\sum_{j=0}^{n-1}\left( 
\begin{array}{c}
n-1 \\ 
j%
\end{array}%
\right) S(j,v-1),
\end{equation*}%
and%
\begin{equation*}
S(n,v)=vS(n-1,v)+S(n-1,v-1),
\end{equation*}%
cf. (\cite{AgohDilcher}, \cite{Carlitz1976}, \cite{Comtet}, \cite%
{LuoSrivatava2010}, \cite{SimsekSpringer}, \cite{YsimsekStirling}).
\end{remark}

The generalized $\lambda $-Stirling type numbers of the second kind can also
be defined by equation (\ref{s5}):

\begin{theorem}
\label{T3}Let $k\in \mathbb{N}_{0}$ and $\lambda \in \mathbb{C}$.%
\begin{equation}
\lambda ^{x}\left( \ln b^{x}\right) ^{m}=\sum_{l=0}^{m}\sum_{j=0}^{\infty
}\left( 
\begin{array}{c}
m \\ 
l%
\end{array}%
\right) \left( 
\begin{array}{c}
x \\ 
j%
\end{array}%
\right) j!\mathcal{S}(l,j;a,b;\lambda )\left( \ln \left( a^{(x-j)}\right)
\right) ^{m-l}.  \label{s5}
\end{equation}
\end{theorem}

\begin{proof}
By using (\ref{s1}), we get%
\begin{equation*}
\left( \lambda b^{t}\right) ^{x}=\sum_{j=0}^{\infty }\left( 
\begin{array}{c}
x \\ 
j%
\end{array}%
\right) j!\sum_{m=0}^{\infty }\mathcal{S}(m,j;a,b;\lambda )\frac{t^{m}}{m!}%
\sum_{n=0}^{\infty }(\ln a^{x-j})^{n}\frac{t^{n}}{n!}.
\end{equation*}%
From the above equation, we obtain%
\begin{equation*}
\lambda ^{x}\sum_{m=0}^{\infty }\left( \ln b\right) ^{m}\frac{t^{m}}{m!}%
=\sum_{m=0}^{\infty }\sum_{j=0}^{\infty }\left( 
\begin{array}{c}
x \\ 
j%
\end{array}%
\right) j!\mathcal{S}(m,j;a,b;\lambda )\frac{t^{m}}{m!}\sum_{n=0}^{\infty
}(\ln a^{x-j})^{n}\frac{t^{n}}{n!}.
\end{equation*}%
Therefore%
\begin{equation*}
\lambda ^{x}\sum_{m=0}^{\infty }\left( \ln b\right) ^{m}\frac{t^{m}}{m!}%
=\sum_{m=0}^{\infty }\left( \sum_{l=0}^{m}\sum_{j=0}^{\infty }\left( 
\begin{array}{c}
m \\ 
l%
\end{array}%
\right) \left( 
\begin{array}{c}
x \\ 
j%
\end{array}%
\right) j!\mathcal{S}(l,j;a,b;\lambda )\left( \ln a^{(x-j)}\right)
^{m-l}\right) \frac{t^{m}}{m!}.
\end{equation*}%
Comparing the coefficients of $\frac{t^{m}}{m!}$ on both sides of the above
equation, we arrive at the desired result.
\end{proof}

\begin{remark}
For $a=0$ and $b=e$, the formula (\ref{s5}) can easily be shown to be
reduced to the following result which is given by Luo and Srivastava \cite[%
Theorem 9]{LuoSrivatava2010}:%
\begin{equation*}
\lambda ^{x}x^{n}=\sum_{l=0}^{\infty }\left( 
\begin{array}{c}
x \\ 
l%
\end{array}%
\right) l!S(n,l;\lambda ),
\end{equation*}%
where $n\in \mathbb{N}_{0}$ and $\lambda \in \mathbb{C}$. For $\lambda =1$,
the above formula is reduced to%
\begin{equation*}
x^{n}=\sum_{v=0}^{n}\left( 
\begin{array}{c}
x \\ 
v%
\end{array}%
\right) v!S(n,v)
\end{equation*}%
cf. (\cite{AgohDilcher}, \cite{Carlitz1976}, \cite{Comtet}, \cite{T. Kim}, 
\cite{LuoSrivatava2010}).
\end{remark}

\section{Generalized array type polynomials}

By using the same motivation with the $\lambda $-Stirling type numbers of
the second kind, we also construct a novel generating function, related to
nonnegative real parameters, of the \textit{generalized array type
polynomials}. We derive some elementary properties including recurrence
relations of these polynomials. The following definition provides a natural
generalization and unification of the array polynomials:

\begin{definition}
Let $a$, $b\in \mathbb{R}^{+}$ ($a\neq b$), $x\in \mathbb{R}$, $\lambda \in 
\mathbb{C}$ and $v\in \mathbb{N}_{0}$. The generalized array type
polynomials $\mathcal{S}_{v}^{n}(x;a,b;\lambda )$\ can be defined by%
\begin{equation}
\mathcal{S}_{v}^{n}(x;a,b;\lambda )=\frac{1}{v!}\sum_{j=0}^{v}(-1)^{v-j}%
\left( 
\begin{array}{c}
v \\ 
j%
\end{array}%
\right) \lambda ^{j}\left( \ln \left( a^{v-j}b^{x+j}\right) \right) ^{n}.
\label{as2}
\end{equation}
\end{definition}

By using the formula (\ref{as2}), we can compute some values of the
polynomials $\mathcal{S}_{v}^{n}(x;a,b;\lambda )$ as follows:%
\begin{equation*}
\mathcal{S}_{0}^{n}(x;a,b;\lambda )=\left( \ln \left( b^{x}\right) \right)
^{n},
\end{equation*}%
\begin{equation*}
\mathcal{S}_{v}^{0}(x;a,b;\lambda )=\frac{\left( 1-\lambda \right) ^{v}}{v!}
\end{equation*}%
and%
\begin{equation*}
\mathcal{S}_{1}^{1}(x;a,b;\lambda )=-\ln (ab^{x})+\lambda \ln (b^{x+1}).
\end{equation*}

\begin{remark}
The polynomials $\mathcal{S}_{v}^{n}(x;a,b;\lambda )$ may be also called
generalized $\lambda $-array type polynomials. By substituting $x=0$ into (%
\ref{as2}), we arrive at (\ref{as1a}):%
\begin{equation*}
\mathcal{S}_{v}^{n}(0;a,b;\lambda )=\mathcal{S}(n,v;a,b;\lambda ).
\end{equation*}%
Setting $a=\lambda =1$ and $b=e$ in (\ref{as2}), we have%
\begin{equation*}
S_{v}^{n}(x)=\frac{1}{v!}\sum_{j=0}^{v}(-1)^{v-j}\left( 
\begin{array}{c}
v \\ 
j%
\end{array}%
\right) \left( x+j\right) ^{n},
\end{equation*}%
a result due to Chang and Ha \cite[Eq-(3.1)]{Chan}, Simsek \cite%
{SimsekSpringer}. It is easy to see that%
\begin{equation*}
S_{0}^{0}(x)=S_{n}^{n}(x)=1,
\end{equation*}%
\begin{equation*}
S_{0}^{n}(x)=x^{n}
\end{equation*}%
and for $v>n$,%
\begin{equation*}
S_{v}^{n}(x)=0
\end{equation*}%
cf. \cite[Eq-(3.1)]{Chan}.
\end{remark}

Generating functions for the polynomial $\mathcal{S}_{v}^{n}(x;a,b,c;\lambda
)$ can be defined as follows:

\begin{definition}
Let $a$, $b\in \mathbb{R}^{+}$ ($a\neq b$), $\lambda \in \mathbb{C}$ and $%
v\in \mathbb{N}_{0}$. The generalized array type polynomials $\mathcal{S}%
_{v}^{n}(x;a,b;\lambda )$\ are defined by means of the following generating
function:%
\begin{equation}
g_{v}(x,t;a,b;\lambda )=\sum_{n=0}^{\infty }\mathcal{S}_{v}^{n}(x;a,b;%
\lambda )\frac{t^{n}}{n!}.  \label{ab1}
\end{equation}
\end{definition}

\begin{theorem}
Let $a$, $b\in \mathbb{R}^{+}$, ($a\neq b$), $\lambda \in \mathbb{C}$ and $%
v\in \mathbb{N}_{0}$.%
\begin{equation}
g_{v}(x,t;a,b;\lambda )=\frac{1}{v!}\left( \lambda b^{t}-a^{t}\right)
^{v}b^{xt}.  \label{ab0}
\end{equation}
\end{theorem}

\begin{proof}
By substituting (\ref{as2}) into the right hand side of (\ref{ab1}), we
obtain%
\begin{equation*}
\sum_{n=0}^{\infty }\mathcal{S}_{v}^{n}(x;a,b;\lambda )\frac{t^{n}}{n!}%
=\sum_{n=0}^{\infty }\left( \frac{1}{v!}\sum_{j=0}^{v}(-1)^{v-j}\left( 
\begin{array}{c}
v \\ 
j%
\end{array}%
\right) \lambda ^{j}\left( \ln \left( a^{v-j}b^{x+j}\right) \right)
^{n}\right) \frac{t^{n}}{n!}.
\end{equation*}%
Therefore%
\begin{equation*}
\sum_{n=0}^{\infty }\mathcal{S}_{v}^{n}(x;a,b;\lambda )\frac{t^{n}}{n!}=%
\frac{1}{v!}\sum_{j=0}^{v}(-1)^{v-j}\left( 
\begin{array}{c}
v \\ 
j%
\end{array}%
\right) \lambda ^{j}\sum_{n=0}^{\infty }\left( \ln \left(
a^{v-j}b^{x+j}\right) \right) ^{n}\frac{t^{n}}{n!}.
\end{equation*}%
The right hand side of the above equation is the Taylor series for $e^{(\ln
\left( a^{v-j}b^{x+j}\right) )t}$, thus we get%
\begin{equation*}
\sum_{n=0}^{\infty }\mathcal{S}_{v}^{n}(x;a,b;\lambda )\frac{t^{n}}{n!}%
=\left( \frac{1}{v!}\sum_{j=0}^{v}(-1)^{v-j}\left( 
\begin{array}{c}
v \\ 
j%
\end{array}%
\right) \lambda ^{j}a^{\left( v-j\right) t}b^{jt}\right) b^{xt}.
\end{equation*}

By using (\ref{s1}) and binomial theorem in the above equation, we arrive at
the desired result.
\end{proof}

\begin{remark}
If we set $\lambda =1$ in (\ref{ab0}), we arrive a new special case of the
array polynomials given by%
\begin{equation*}
f_{S,v}(t;a,b)b^{tx}=\sum_{n=0}^{\infty }\mathcal{S}_{v}^{n}(x;a,b)\frac{%
t^{n}}{n!}.
\end{equation*}%
In the special case when%
\begin{equation*}
a=\lambda =1\text{ and }b=e,
\end{equation*}%
the generalized array polynomials $\mathcal{S}_{v}^{n}(x;a,b;\lambda )$
defined by (\ref{ab0}) would lead us at once to the classical array
polynomials $\mathcal{S}_{v}^{n}(x)$, which are defined by means of the
following generating function:%
\begin{equation*}
\frac{\left( e^{t}-1\right) ^{v}}{v!}e^{tx}=\sum_{n=0}^{\infty }S_{v}^{n}(x)%
\frac{t^{n}}{n!},
\end{equation*}%
which yields to the generating function for the array polynomials $%
S_{v}^{n}(x)$ studied by Chang and Ha \cite{Chan} see also cf. (\cite{cagic}%
, \cite{SimsekSpringer}).
\end{remark}

The polynomials $\mathcal{S}_{v}^{n}(x;a,b;\lambda )$\ defined by (\ref{ab0}%
) have many interesting properties which we give in this section.

We set%
\begin{equation}
g_{v}(x,t;a,b;\lambda )=b^{xt}f_{S,v}(t;a,b;\lambda ).  \label{1Sse}
\end{equation}

\begin{theorem}
The following formula holds true:%
\begin{equation}
\mathcal{S}_{v}^{n}(x;a,b;\lambda )=\sum_{j=0}^{n}\left( 
\begin{array}{c}
n \\ 
j%
\end{array}%
\right) \mathcal{S}(j,v;a,b;\lambda )\left( \ln b^{x}\right) ^{n-j}.
\label{1Ssc}
\end{equation}
\end{theorem}

\begin{proof}
By using (\ref{1Sse}), we obtain%
\begin{equation*}
\sum_{n=0}^{\infty }\mathcal{S}_{v}^{n}(x;a,b;\lambda )\frac{t^{n}}{n!}%
=\sum_{n=0}^{\infty }\mathcal{S}(n,v;a,b;\lambda )\frac{t^{n}}{n!}%
\sum_{n=0}^{\infty }\left( \ln b^{x}\right) ^{n}\frac{t^{n}}{n!}.
\end{equation*}

From the above equation, we get%
\begin{equation*}
\sum_{n=0}^{\infty }\mathcal{S}_{v}^{n}(x;a,b;\lambda )\frac{t^{n}}{n!}%
=\sum_{n=0}^{\infty }\left( \sum_{j=0}^{n}\left( 
\begin{array}{c}
n \\ 
j%
\end{array}%
\right) \mathcal{S}(j,v;a,b)\left( \ln b^{x}\right) ^{n-j}\right) \frac{t^{n}%
}{n!}.
\end{equation*}

Comparing the coefficients of $t^{n}$ on both sides of the above equation,
we arrive at the desired result.
\end{proof}

\begin{remark}
In the special case when $a=\lambda =1$ and $b=e$, equation (\ref{1Ssc}) is
reduced to%
\begin{equation*}
S_{v}^{n}(x)=\sum_{j=0}^{n}\left( 
\begin{array}{c}
n \\ 
j%
\end{array}%
\right) x^{n-j}S(j,v)
\end{equation*}%
cf. \cite[Theorem 2]{SimsekSpringer}.
\end{remark}

By differentiating $j$ times both sides of (\ref{ab0}) with respect to the
variable $x$, we obtain the following differential equation:%
\begin{equation*}
\frac{\partial ^{j}}{\partial x^{j}}g_{v}(x,t;a,b;\lambda )=t^{j}\left( \ln
b\right) ^{j}g_{v}(x,t;a,b;\lambda ).
\end{equation*}

From this equation, we arrive at higher order derivative of \ the array type
polynomials by the following theorem:

\begin{theorem}
\label{TEo2} Let $n$, $j\in \mathbb{N}$ with $j\leq n$. Then we have%
\begin{equation*}
\frac{\partial ^{j}}{\partial x^{j}}\mathcal{S}_{v}^{n}(x;a,b;\lambda
)=\left\{ n\right\} _{j}\left( \ln (b)\right) ^{j}\mathcal{S}%
_{v}^{n-j}(x;a,b;\lambda ).
\end{equation*}
\end{theorem}

\begin{remark}
By setting $a=\lambda =j=1$ and $b=e$ in Theorem \ref{TEo2}, we have%
\begin{equation*}
\frac{d}{dx}S_{v}^{n}(x)=nS_{v}^{n-1}(x)
\end{equation*}%
cf. \cite{SimsekSpringer}.
\end{remark}

From (\ref{ab0}), we get the following functional equation:%
\begin{equation}
g_{v_{1}}(x_{1},t;a,b;\lambda )g_{v_{2}}(x_{2},t;a,b;\lambda )=\left( 
\begin{array}{c}
v_{1}+v_{2} \\ 
v_{1}%
\end{array}%
\right) g_{v_{1}+v_{2}}(x_{1}+x_{2},t;a,b;\lambda ).  \label{as3}
\end{equation}%
From this functional equation, we obtain the following identity:

\begin{theorem}
\begin{equation*}
\left( 
\begin{array}{c}
v_{1}+v_{2} \\ 
v_{1}%
\end{array}%
\right) \mathcal{S}_{v_{1}+v_{2}}^{n}(x_{1}+x_{2};a,b;\lambda
)=\dsum\limits_{j=0}^{n}\left( 
\begin{array}{c}
n \\ 
j%
\end{array}%
\right) \mathcal{S}_{v_{1}}^{j}(x_{1};a,b;\lambda )\mathcal{S}%
_{v_{2}}^{n-j}(x_{2};a,b;\lambda ).
\end{equation*}
\end{theorem}

\begin{proof}
Combining (\ref{ab1}) and (\ref{as3}), we get%
\begin{eqnarray*}
&&\sum_{n=0}^{\infty }\mathcal{S}_{v_{1}}^{n}(x_{1};a,b;\lambda )\frac{t^{n}%
}{n!}\sum_{n=0}^{\infty }\mathcal{S}_{v_{2}}^{n}(x_{2};a,b;\lambda )\frac{%
t^{n}}{n!} \\
&=&\left( 
\begin{array}{c}
v_{1}+v_{2} \\ 
v_{1}%
\end{array}%
\right) \sum_{n=0}^{\infty }\mathcal{S}_{v_{1}+v_{2}}^{n}(x_{1}+x_{2};a,b;%
\lambda )\frac{t^{n}}{n!}.
\end{eqnarray*}%
Therefore%
\begin{eqnarray*}
&&\sum_{n=0}^{\infty }\left( \dsum\limits_{j=0}^{n}\left( 
\begin{array}{c}
n \\ 
j%
\end{array}%
\right) \mathcal{S}_{v_{1}}^{j}(x_{1};a,b;\lambda )\mathcal{S}%
_{v_{2}}^{n-j}(x_{2};a,b;\lambda )\right) \frac{t^{n}}{n!} \\
&=&\left( 
\begin{array}{c}
v_{1}+v_{2} \\ 
v_{1}%
\end{array}%
\right) \sum_{n=0}^{\infty }\mathcal{S}_{v_{1}+v_{2}}^{n}(x_{1}+x_{2};a,b;%
\lambda )\frac{t^{n}}{n!}\text{.}
\end{eqnarray*}%
Comparing the coefficients of $\frac{t^{n}}{n!}$ on both sides of the above
equation, we arrive at the desired result.
\end{proof}

\section{Generalized Eulerian type numbers and polynomials}

In this section, we provide generating functions, related to nonnegative
real parameters, for the generalized Eulerian type polynomials and numbers,
that is, the so called \textit{generalized Apostol type Frobenius Euler
polynomials} \textit{and numbers}. We derive fundamental properties,
recurrence relations and many new identities for these polynomials and
numbers based on the generating functions, functional equations and
differential equations.

These polynomials and numbers have many applications in many branches of
Mathematics.

The following definition gives us a natural generalization of the Eulerian
polynomials:

\begin{definition}
Let $a,$ $b\in \mathbb{R}^{+}$ $(a\neq b),$ $x\in \mathbb{R},$ $\lambda \in 
\mathbb{C}$ and $u\in \mathbb{C\diagdown }\left\{ 1\right\} $. The\
generalized Eulerian type polynomials $\mathcal{H}_{n}(x;u;a,b,c;\lambda )$
are defined by means of the following generating function:%
\begin{equation}
F_{\lambda }(t,x;u,a,b,c)=\frac{\left( a^{t}-u\right) c^{xt}}{\lambda b^{t}-u%
}=\sum_{n=0}^{\infty }\mathcal{H}_{n}(x;u;a,b,c;\lambda )\frac{t^{n}}{n!}.
\label{4ge1}
\end{equation}
\end{definition}

By substituting $x=0$ into (\ref{4ge1}), we obtain%
\begin{equation*}
\mathcal{H}_{n}(0;u;a,b,c;\lambda )=\mathcal{H}_{n}(u;a,b,c;\lambda ),
\end{equation*}%
where $\mathcal{H}_{n}(u;a,b,c;\lambda )$ denotes \textit{generalized
Eulerian type numbers}.

\begin{remark}
Substituting $a=1$ into (\ref{4ge1}), we have%
\begin{equation*}
\frac{\left( 1-u\right) c^{xt}}{\lambda b^{t}-u}=\sum_{n=0}^{\infty }%
\mathcal{H}_{n}(x;u;1,b,c;\lambda )\frac{t^{n}}{n!}
\end{equation*}%
a result due to Kurt and Simsek \cite{burakSimsek}. In their special case
when $\lambda =1$ and $b=c=e$, the \textit{generalized }Eulerian type
polynomials $\mathcal{H}_{n}(x;u;1,b,c;\lambda )$\ are reduced to the
Eulerian polynomials or Frobenius Euler polynomials which are defined by
means of the following generating function:%
\begin{equation}
\frac{\left( 1-u\right) e^{xt}}{e^{t}-u}=\sum_{n=0}^{\infty }H_{n}(x;u)\frac{%
t^{n}}{n!},  \label{mt2}
\end{equation}%
with, of course, $H_{n}(0;u)=H_{n}(u)$ denotes the so-called Eulerian
numbers cf. (\cite{Carlitz}, \cite{Carlitz1952}, \cite{Carlitz1953G}, \cite%
{Carlitz1976}, \cite{KimmskimlcjangJIA}, \cite{YsimsekKim}, \cite%
{KimSimskJKM}, \cite{SimsekBKMS}, \cite{SimsekJNT}, \cite{srivas18}, \cite%
{Tsumura}). Substituting $u=-1$, into (\ref{mt2}), we have%
\begin{equation*}
H_{n}(x;-1)=E_{n}(x)
\end{equation*}%
where $E_{n}(x)$\ denotes Euler polynomials which are defined by means of
the following generating function:%
\begin{equation}
\frac{2e^{xt}}{e^{t}+1}=\sum_{n=0}^{\infty }E_{n}(x)\frac{t^{n}}{n!}
\label{1Ssf}
\end{equation}%
where $\left\vert t\right\vert <\pi $\ cf. \cite{AgohDilcher}-\cite{walum}.
\end{remark}

The following elementary properties of the generalized Eulerian type
polynomials and numbers are derived from their generating functions in (\ref%
{4ge1}).

\begin{theorem}
(\textit{Recurrence relation} for the generalized Eulerian type numbers):
For $n=0$, we have%
\begin{equation*}
\mathcal{H}_{0}(u;a,b;\lambda )=\left\{ 
\begin{array}{c}
\frac{1-u}{\lambda -u}\text{ if }a=1, \\ 
\\ 
\frac{u}{\lambda -u}\text{ if }a\neq 1.%
\end{array}%
\right.
\end{equation*}%
For $n>0$, following the usual convention of symbolically replacing $\left( 
\mathcal{H}(u;a,b;\lambda )\right) ^{n}$ by $\mathcal{H}_{n}(u;a,b;\lambda )$%
, we have%
\begin{equation*}
\lambda \left( \ln b+\mathcal{H}(u;a,b;\lambda )\right) ^{n}-u\mathcal{H}%
_{n}(u;a,b;\lambda )=\left( \ln a\right) ^{n}.
\end{equation*}
\end{theorem}

\begin{proof}
By using (\ref{4ge1}), we obtain%
\begin{equation*}
\sum_{n=0}^{\infty }\left( \ln a\right) ^{n}\frac{t^{n}}{n!}%
-u=\sum_{n=0}^{\infty }\left( \lambda \left( \ln b+\mathcal{H}(u;a,b;\lambda
)\right) ^{n}-u\mathcal{H}_{n}(u;a,b;\lambda )\right) \frac{t^{n}}{n!}.
\end{equation*}%
Comparing the coefficients of $\frac{t^{n}}{n!}$ on both sides of the above
equation, we arrive at the desired result.
\end{proof}

By differentiating both sides of equation (\ref{4ge1}) with respect to the
variable $x$, we obtain the following higher order differential equation:%
\begin{equation}
\frac{\partial ^{j}}{\partial x^{j}}F_{\lambda }(t,x;u,a,b,c)=\left( \ln
\left( c^{t}\right) \right) ^{j}F_{\lambda }(t,x;u,a,b,c).  \label{F1}
\end{equation}%
From this equation, we arrive at higher order derivative of \ the
generalized Eulerian type polynomials by the following theorem:

\begin{theorem}
\label{TEo3} Let $n$, $j\in \mathbb{N}$ with $j\leq n$. Then we have%
\begin{equation*}
\frac{\partial ^{j}}{\partial x^{j}}\mathcal{H}_{n}(x;u;a,b,c;\lambda
)=\left\{ n\right\} _{j}\left( \ln \left( c\right) \right) ^{j}\mathcal{H}%
_{n-j}(x;u;a,b,c;\lambda ).
\end{equation*}
\end{theorem}

\begin{proof}
Combining (\ref{4ge1}) and (\ref{F1}), we have%
\begin{equation*}
\sum_{n=0}^{\infty }\frac{\partial ^{j}}{\partial x^{j}}\mathcal{H}%
_{n}(x;u;a,b,c;\lambda )\frac{t^{n}}{n!}=\left( \ln c\right)
^{j}\sum_{n=0}^{\infty }\mathcal{H}_{n}(x;u;a,b,c;\lambda )\frac{t^{n+j}}{n!}%
.
\end{equation*}%
From the above equation, we get%
\begin{equation*}
\sum_{n=0}^{\infty }\frac{\partial ^{j}}{\partial x^{j}}\mathcal{H}%
_{n}(x;u;a,b,c;\lambda )\frac{t^{n}}{n!}=\left( \ln c\right)
^{j}\sum_{n=0}^{\infty }\left\{ n\right\} _{j}\mathcal{H}_{n-j}(x;u;a,b,c;%
\lambda )\frac{t^{n}}{n!}.
\end{equation*}%
Comparing the coefficients of $\frac{t^{n}}{n!}$ on both sides of the above
equation, we arrive at the desired result.
\end{proof}

\begin{remark}
Setting $j=1$ in Theorem \ref{TEo3}, we have%
\begin{equation*}
\frac{\partial }{\partial x}\mathcal{H}_{n}(x;u;a,b,c;\lambda )=n\mathcal{H}%
_{n-1}(x;u;a,b,c;\lambda )\ln \left( c\right) .
\end{equation*}%
In their special case when%
\begin{equation*}
a=\lambda =1\text{ and }b=c=e,
\end{equation*}%
Theorem \ref{TEo3}\ is reduced to the following well known result:%
\begin{equation*}
\frac{\partial ^{j}}{\partial x^{j}}H_{n}(x;u)=\frac{n!}{(n-j)!}H_{n-j}(x;u)
\end{equation*}%
cf. \cite[Eq-(3.5)]{Carlitz}. Substituting $j=1$ into the above equation, we
have%
\begin{equation*}
\frac{\partial }{\partial x}H_{n}(x;u)=nH_{n-1}(x;u)
\end{equation*}%
cf. (\cite[Eq-(3.5)]{Carlitz}, \cite{burakSimsek}).
\end{remark}

\begin{theorem}
\label{t8} The following explicit representation formula holds true:%
\begin{eqnarray*}
&&\left( x\ln c+\ln a\right) ^{n}-ux^{n}\left( \ln c\right) ^{n} \\
&=&\lambda \left( x\ln c+\ln b+\mathcal{H}(u;a,b;\lambda )\right)
^{n}-u\left( x\ln c+\mathcal{H}(u;a,b;\lambda )\right) ^{n}.
\end{eqnarray*}
\end{theorem}

\begin{proof}
By using (\ref{4ge1}) and the \textit{umbral calculus convention}, we obtain%
\begin{equation*}
\frac{a^{t}-u}{\lambda b^{t}-u}=e^{H\left( u;a,b;\lambda \right) t}.
\end{equation*}%
From the above equation, we get%
\begin{eqnarray*}
&&\sum_{n=0}^{\infty }\left( \left( \ln a+x\ln c\right) ^{n}-u\left( x\ln
c\right) \right) \frac{t^{n}}{n!} \\
&=&\sum_{n=0}^{\infty }\left( \lambda \left( \mathcal{H}\left( u;a,b;\lambda
\right) +\ln b+x\ln c\right) ^{n}-u\left( \mathcal{H}_{n}\left(
u;a,b;\lambda \right) +x\ln c\right) ^{n}\right) \frac{t^{n}}{n!}.
\end{eqnarray*}%
Comparing the coefficients of $\frac{t^{n}}{n!}$ on both sides of the above
equation, we arrive at the desired result.
\end{proof}

\begin{remark}
By substituting $a=\lambda =1$ and $b=c=e$ into Theorem \ref{t8}, we have%
\begin{equation}
\left( 1-u\right) x^{n}=H_{n}(x+1;u)-uH_{n}(x;u)  \label{at8}
\end{equation}%
cf. (\cite[Eq-(3.3)]{Carlitz}, \cite{Tsumura}). By setting $u=-1$ in the
above equation, we have%
\begin{equation*}
2x^{n}=E_{n}(x+1)+E_{n}(x)
\end{equation*}%
a result due to Shiratani \cite{K. Shiratani}. By using (\ref{at8}), Carlitz 
\cite{Carlitz} studied on the \textbf{Mirimonoff polynomial} $f_{n}(0,m)$
which is defined by%
\begin{eqnarray*}
f_{n}(x,m) &=&\dsum\limits_{j=0}^{m-1}(x+j)^{n}u^{m-j-1} \\
&=&\frac{H_{n}(x+m;u)-u^{m}H_{n}(x;u)}{1-u}.
\end{eqnarray*}%
By applying Theorem \ref{t8}, one may generalize the Mirimonoff polynomial.
\end{remark}

\begin{theorem}
The following explicit representation formula holds true:%
\begin{equation}
\mathcal{H}_{n}(x;u;a,b,c;\lambda )=\sum_{j=0}^{n}\left( 
\begin{array}{c}
n \\ 
j%
\end{array}%
\right) \left( x\ln c\right) ^{n-j}\mathcal{H}_{j}(u;a,b,c;\lambda ).
\label{Te}
\end{equation}
\end{theorem}

\begin{proof}
By using (\ref{4ge1}), we get%
\begin{equation*}
\sum_{n=0}^{\infty }\mathcal{H}_{n}(u;a,b,c;\lambda )\frac{t^{n}}{n!}%
\sum_{n=0}^{\infty }\left( \ln c\right) ^{n}\frac{t^{n}}{n!}%
=\sum_{n=0}^{\infty }\mathcal{H}_{n}(x;u;a,b,c;\lambda )\frac{t^{n}}{n!}.
\end{equation*}%
From the above equation, we obtain%
\begin{equation*}
\sum_{n=0}^{\infty }\left( \sum_{j=0}^{n}\left( 
\begin{array}{c}
n \\ 
j%
\end{array}%
\right) \left( x\ln c\right) ^{n-j}\mathcal{H}_{j}(u;a,b,c;\lambda )\right) 
\frac{t^{n}}{n!}=\sum_{n=0}^{\infty }\mathcal{H}_{n}(x;u;a,b,c;\lambda )%
\frac{t^{n}}{n!}.
\end{equation*}%
Comparing the coefficients of $\frac{t^{n}}{n!}$ on both sides of the above
equation, we arrive at the desired result.
\end{proof}

\begin{remark}
Substituting $a=\lambda =1$ and $b=c=e$ into (\ref{Te}), we have%
\begin{equation*}
H_{n}(x;u)=\sum_{j=0}^{n}\left( 
\begin{array}{c}
n \\ 
j%
\end{array}%
\right) x^{n-j}H_{j}(u)
\end{equation*}%
cf. (\cite{Carlitz}, \cite{Carlitz1952}, \cite{Carlitz1953G}, \cite%
{Carlitz1976}, \cite{KimmskimlcjangJIA}, \cite{YsimsekKim}, \cite%
{KimSimskJKM}, \cite{burakSimsek}, \cite{SimsekBKMS}, \cite{SimsekJNT}, \cite%
{srivas18}, \cite{Tsumura}).
\end{remark}

\begin{remark}
From (\ref{Te}), we easily get%
\begin{equation*}
\mathcal{H}_{n}(x;u;a,b,c;\lambda )=\left( \mathcal{H}(u;a,b,c;\lambda
)+x\ln c\right) ^{n},
\end{equation*}%
where after expansion of the right member, $\mathcal{H}^{n}(u;a,b,c;\lambda
) $ is replaced by $\mathcal{H}_{n}(u;a,b,c;\lambda )$, we use this
convention frequently throughout of this paper.
\end{remark}

\begin{theorem}
\begin{equation}
\mathcal{H}_{n}(x+y;u;a,b,c;\lambda )=\sum_{j=0}^{n}\left( 
\begin{array}{c}
n \\ 
j%
\end{array}%
\right) \left( y\ln c\right) ^{n-j}\mathcal{H}_{j}(x;u;a,b,c;\lambda ).
\label{1Ssa}
\end{equation}
\end{theorem}

\begin{proof}
By using (\ref{4ge1}), we have%
\begin{equation*}
\sum_{n=0}^{\infty }\mathcal{H}_{n}(x+y;u;a,b,c;\lambda )\frac{t^{n}}{n!}%
=\sum_{n=0}^{\infty }\left( y\ln c\right) ^{n}\frac{t^{n}}{n!}%
.\sum_{n=0}^{\infty }\mathcal{H}_{n}(x;u;a,b,c;\lambda )\frac{t^{n}}{n!}.
\end{equation*}%
Therefore%
\begin{equation*}
\sum_{n=0}^{\infty }\mathcal{H}_{n}(x+y;u;a,b,c;\lambda )\frac{t^{n}}{n!}%
=\sum_{n=0}^{\infty }\sum_{j=0}^{n}\left( 
\begin{array}{c}
n \\ 
j%
\end{array}%
\right) \left( y\ln c\right) ^{n-j}\mathcal{H}_{j}(x,u;a,b,c;\lambda )\frac{%
t^{n}}{n!}.
\end{equation*}%
Comparing the coefficients of $\frac{t^{n}}{n!}$ on both sides of the above
equation, we arrive at the desired result.
\end{proof}

\begin{remark}
In the special case when $a=\lambda =1$ and $b=c=e$, equation (\ref{1Ssa})\
is reduced to the following result:%
\begin{equation*}
H_{n}(x+y)=\sum_{j=0}^{n}\left( 
\begin{array}{c}
n \\ 
j%
\end{array}%
\right) y^{n-j}H_{j}(x,u)
\end{equation*}%
cf. \cite[Eq-(3.6)]{Carlitz}. Substituting $u=-1$ into the above equation,
we get the following well-known result:%
\begin{equation}
E_{n}(x+y)=\sum_{j=0}^{n}\left( 
\begin{array}{c}
n \\ 
j%
\end{array}%
\right) y^{n-j}E_{j}(x).  \label{4Eq}
\end{equation}
\end{remark}

By using (\ref{4ge1}), we define the following functional equation:%
\begin{equation}
F_{\lambda ^{2}}(t,x;u^{2},a^{2},b^{2},c)c^{yt}=F_{\lambda
}(t,x;u,a,b,c)F_{\lambda }(t,y;-u,a,b,c).  \label{4EqY}
\end{equation}

\begin{theorem}
\begin{equation}
\mathcal{H}_{n}(x+y;u^{2};a,b,c;\lambda ^{2})=\left( \mathcal{H}%
(x;u;a,b,c;\lambda )+\mathcal{H}(y;-u;a,b,c;\lambda )\right) ^{n}.
\label{4Eqy1}
\end{equation}
\end{theorem}

\begin{proof}
Combining (\ref{4EqY}) and (\ref{1Ssa}), we easily arrive at the desired
result.
\end{proof}

\begin{remark}
In the special case when $a=\lambda =1$ and $b=c=e$, equation (\ref{4Eqy1})\
is reduced to the following result:%
\begin{equation*}
H_{n}(x+y;u^{2})=\sum_{j=0}^{n}\left( 
\begin{array}{c}
n \\ 
j%
\end{array}%
\right) H_{j}(x;u)H_{n-j}(y;-u)
\end{equation*}%
cf. \cite[Eq-(3.17)]{Carlitz}.
\end{remark}

\begin{theorem}
\label{TeoE}%
\begin{equation*}
(-1)^{n}\mathcal{H}_{n}(1-x;u^{-1};a,b,c;\lambda ^{-1})=\lambda
\sum_{j=0}^{n}\left( 
\begin{array}{c}
n \\ 
j%
\end{array}%
\right) \left( \ln \left( \frac{b}{a}\right) \right) ^{n-j}\mathcal{H}%
_{j}(x-1,u;a,b,c;\lambda ).
\end{equation*}
\end{theorem}

\begin{proof}
By using (\ref{4ge1}), we obtain%
\begin{equation*}
\frac{\left( a^{-t}-u^{-1}\right) c^{-\left( 1-x\right) t}}{\lambda
^{-1}b^{-t}-u^{-1}}=\lambda \left( \frac{b}{a}\right) ^{t}\sum_{n=0}^{\infty
}\mathcal{H}_{n}(x-1;u;a,b,c;\lambda )\frac{t^{n}}{n!}.
\end{equation*}%
From the above equation, we get%
\begin{eqnarray*}
&&\sum_{n=0}^{\infty }\mathcal{H}_{n}(1-x;u^{-1};a,b,c;\lambda ^{-1})\frac{%
(-1)^{n}t^{n}}{n!} \\
&=&\lambda \left( \sum_{n=0}^{\infty }\mathcal{H}_{n}(x-1;u;a,b,c;\lambda )%
\frac{t^{n}}{n!}\right) \left( \sum_{n=0}^{\infty }\left( \ln \left( \frac{b%
}{a}\right) \right) ^{n}\frac{t^{n}}{n!}\right) .
\end{eqnarray*}%
Therefore%
\begin{eqnarray*}
&&\sum_{n=0}^{\infty }(-1)^{n}\mathcal{H}_{n}(1-x;u^{-1};a,b,c;\lambda ^{-1})%
\frac{t^{n}}{n!} \\
&=&\sum_{n=0}^{\infty }\left( \lambda \sum_{j=0}^{n}\left( 
\begin{array}{c}
n \\ 
j%
\end{array}%
\right) \left( \ln \left( \frac{b}{a}\right) \right) ^{n-j}\mathcal{H}%
_{j}(x-1,u;a,b,c;\lambda )\right) \frac{t^{n}}{n!}.
\end{eqnarray*}%
Comparing the coefficients of $\frac{t^{n}}{n!}$ on both sides of the above
equation, we arrive at the desired result.
\end{proof}

\begin{remark}
In their special case when $a=\lambda =1$ and $b=c=e$, Theorem \ref{TeoE}\
is reduced to the following result:%
\begin{equation*}
(-1)^{n}H_{n}(1-x;u^{-1})=H_{n}(x-1,u)
\end{equation*}%
cf. \cite[Eq-(3.7)]{Carlitz}. Substituting $u=-1$ into the above equation,
we get the following well-known result:%
\begin{equation*}
(-1)^{n}E_{n}(1-x)=E_{n}(x)
\end{equation*}%
cf. (\cite[Eq-(3.7)]{Carlitz}, \cite{http}, \cite{Parashar}, \cite{K.
Shiratani}, \cite{Srivastava2011}).
\end{remark}

\begin{theorem}
\label{TeoE-1}%
\begin{equation*}
\mathcal{H}_{n}\left( \frac{x+y}{2};u^{2};a,b,c;\lambda ^{2}\right)
=\sum_{j=0}^{n}\left( 
\begin{array}{c}
n \\ 
j%
\end{array}%
\right) \frac{\mathcal{H}_{j}(x;u;a,b,c;\lambda )\mathcal{H}%
_{n-j}(y;-u;a,b,c;\lambda )}{2^{n}}.
\end{equation*}
\end{theorem}

\begin{proof}
By using (\ref{4ge1}), we get%
\begin{eqnarray*}
&&\sum_{n=0}^{\infty }\mathcal{H}_{n}\left( \frac{x+y}{2};u^{2};a,b,c;%
\lambda ^{2}\right) \frac{2^{n}t^{n}}{n!} \\
&=&\sum_{n=0}^{\infty }\left( \sum_{j=0}^{n}\left( 
\begin{array}{c}
n \\ 
j%
\end{array}%
\right) \mathcal{H}_{j}(x;u;a,b,c;\lambda )\mathcal{H}_{n-j}(y;-u;a,b,c;%
\right) \frac{t^{n}}{n!}.
\end{eqnarray*}%
Comparing the coefficients of $\frac{t^{n}}{n!}$ on both sides of the above
equation, we arrive at the desired result.
\end{proof}

\begin{remark}
When $a=\lambda =1$ and $b=c=e$, Theorem \ref{TeoE-1}\ is reduced to the
following result:%
\begin{equation*}
\mathcal{H}_{n}\left( \frac{x+y}{2};u^{2}\right) =2^{-n}\sum_{j=0}^{n}\left( 
\begin{array}{c}
n \\ 
j%
\end{array}%
\right) \mathcal{H}_{j}(x;u)\mathcal{H}_{n-j}(y;-u),
\end{equation*}%
cf. \cite[Eq-(3.17)]{Carlitz}.
\end{remark}

\subsection{Multiplication formulas for normalized polynomials}

In this section, using generating functions, we derive \textit{%
multiplication formulas} in terms of the normalized polynomials which are
related to the generalized Eulerian type polynomials, the Bernoulli and the
Euler polynomials.

\begin{theorem}
\label{T11}(Multiplication formula) Let $y\in \mathbb{N}$. Then we have%
\begin{eqnarray}
&&\mathcal{H}_{n}(yx;u;a,b,b;\lambda )  \label{mMF} \\
&=&y^{n}\dsum\limits_{k=0}^{n}\dsum\limits_{j=0}^{y-1}\left( 
\begin{array}{c}
n \\ 
k%
\end{array}%
\right) \frac{\lambda ^{j}\left( \ln a\right) ^{n-k}}{u^{j+1-y}-u^{j+1}}%
\mathcal{H}_{k}\left( x+\frac{j}{y};u^{y};a,b,b;\lambda ^{y}\right)  \notag
\\
&&\times \left( H_{n-k}\left( \frac{1}{y};u^{y}\right) -uH_{n-k}\left(
u^{y}\right) \right) ,  \notag
\end{eqnarray}%
where $H_{n}\left( x;u\right) $ and $H_{n}\left( u\right) $ denote the
Eulerian polynomials and numbers, respectively.
\end{theorem}

\begin{proof}
Substituting $c=b$ into (\ref{4ge1}), we have%
\begin{equation}
\dsum\limits_{n=0}^{\infty }\mathcal{H}_{n}(x;u;a,b,b;\lambda )\frac{t^{n}}{%
n!}=\frac{\left( a^{t}-u\right) b^{xt}}{\lambda b^{t}-u}=\left( \frac{a^{t}-u%
}{-u}\right) \frac{\left( a^{t}-u\right) b^{xt}}{1-\frac{\lambda b^{t}}{u}}.
\label{MT1}
\end{equation}%
By using the following finite geometric series%
\begin{equation*}
\dsum\limits_{j=0}^{y-1}\left( \frac{\lambda b^{t}}{u}\right) ^{j}=\frac{%
1-\left( \frac{\lambda b^{t}}{u}\right) ^{y}}{1-\frac{\lambda b^{t}}{u}},
\end{equation*}%
on the right-hand side of (\ref{MT1}), we obtain%
\begin{equation*}
\dsum\limits_{n=0}^{\infty }\mathcal{H}_{n}(x;u;a,b,b;\lambda )\frac{t^{n}}{%
n!}=\frac{\left( a^{t}-u\right) b^{xt}}{-u\left( 1-\left( \frac{\lambda b^{t}%
}{u}\right) ^{y}\right) }\dsum\limits_{j=0}^{y-1}\left( \frac{\lambda b^{t}}{%
u}\right) ^{j}.
\end{equation*}%
From this equation, we get%
\begin{equation*}
\dsum\limits_{n=0}^{\infty }\mathcal{H}_{n}(x;u;a,b,b;\lambda )\frac{t^{n}}{%
n!}=\frac{\left( a^{t}-u\right) }{\left( a^{yt}-u^{y}\right) }%
\dsum\limits_{j=0}^{y-1}\frac{\lambda ^{j}}{u^{j+1-y}}\frac{\left(
a^{yt}-u^{y}\right) b^{yt\left( \frac{x+j}{y}\right) }}{\left( \lambda
b^{yt}-u^{y}\right) }.
\end{equation*}%
Now by making use of the generating functions (\ref{4ge1}) and (\ref{mt2})
on the right-hand side of the above equation, we obtain%
\begin{eqnarray*}
&&\dsum\limits_{n=0}^{\infty }\mathcal{H}_{n}(x;u;a,b,b;\lambda )\frac{t^{n}%
}{n!} \\
&=&\frac{1}{1-u^{y}}\dsum\limits_{j=0}^{y-1}\frac{\lambda ^{j}}{u^{j+1-y}}%
\left( \dsum\limits_{n=0}^{\infty }\mathcal{H}_{n}\left( \frac{x+j}{y}%
;u^{y};a,b,b;\lambda ^{y}\right) \frac{y^{n}t^{n}}{n!}\right) \\
&&\times \left( \dsum\limits_{n=0}^{\infty }\left( H_{n}\left( \frac{1}{y}%
;u^{y}\right) -uH_{n}\left( u^{y}\right) \right) \frac{\left( y\ln a\right)
^{n}t^{n}}{n!}\right) .
\end{eqnarray*}%
Therefore%
\begin{eqnarray*}
&&\dsum\limits_{n=0}^{\infty }\mathcal{H}_{n}(x;u;a,b,b;\lambda )\frac{t^{n}%
}{n!} \\
&=&\dsum\limits_{n=0}^{\infty
}\dsum\limits_{k=0}^{n}\dsum\limits_{j=0}^{y-1}\left( 
\begin{array}{c}
n \\ 
k%
\end{array}%
\right) \frac{y^{n}\lambda ^{j}\left( \ln a\right) ^{n-k}}{u^{j+1-y}-u^{j+1}}%
\mathcal{H}_{k}\left( \frac{x+j}{y};u^{y};a,b,b;\lambda ^{y}\right) \\
&&\times \left( H_{n-k}\left( \frac{1}{y};u^{y}\right) -uH_{n-k}\left(
u^{y}\right) \right) \frac{t^{n}}{n!}.
\end{eqnarray*}%
By equating the coefficients of $\frac{t^{n}}{n!}$\ on both sides, we get%
\begin{eqnarray*}
&&\mathcal{H}_{n}(x;u;a,b,b;\lambda ) \\
&=&\dsum\limits_{k=0}^{n}\dsum\limits_{j=0}^{y-1}\left( 
\begin{array}{c}
n \\ 
k%
\end{array}%
\right) \frac{y^{n}\lambda ^{j}\left( \ln a\right) ^{n-k}}{u^{j+1-y}-u^{j+1}}%
\mathcal{H}_{k}\left( \frac{x+j}{y};u^{y};a,b,b;\lambda ^{y}\right) \\
&&\times \left( H_{n-k}\left( \frac{1}{y};u^{y}\right) -uH_{n-k}\left(
u^{y}\right) \right) .
\end{eqnarray*}%
Finally, by replacing $x$ by $yx$ on both sides of the above equation, we
arrive at the desired result.
\end{proof}

\begin{remark}
By substituting $a=1$ into Theorem \ref{T11}, for $n=k$, we obtain%
\begin{equation}
\mathcal{H}_{n}(yx;u;1,b,b;\lambda )=y^{n}u^{y-1}\frac{1-u}{1-u^{y}}%
\dsum\limits_{j=0}^{y-1}\frac{\lambda ^{j}}{u^{j}}\mathcal{H}_{n}\left( x+%
\frac{j}{y};u^{y};1,b,b;\lambda ^{y}\right) .  \label{mmF}
\end{equation}%
By substituting $b=e$ and $\lambda =1$ into the above equation, we arrive at
the multiplication formula for the Eulerian polynomials%
\begin{equation}
H_{n}(yx;u)=y^{n}u^{y-1}\frac{\left( 1-u\right) }{1-u^{y}}%
\dsum\limits_{j=0}^{y-1}\frac{1}{u^{j}}H_{n}\left( x+\frac{j}{y}%
;u^{y}\right) ,  \label{MF-0}
\end{equation}%
cf. (\cite{Carlitz1953G}, \cite[Eq-(3.12)]{Carlitz}). If $u=-1$, then the
above equation reduces to the well known multiplication formula for the
Euler polynomials: for $y$ is an odd positive integer, we have%
\begin{equation}
E_{n}(yx)=y^{n}\dsum\limits_{j=0}^{y-1}(-1)^{j}E_{n}\left( x+\frac{j}{y}%
\right) ,  \label{MF}
\end{equation}%
where $E_{n}(x)$ denotes the Euler polynomials in the usual notation. If$\ y$
is an even positive integer, we have%
\begin{equation}
E_{n}(yx)=\frac{2y^{n-1}}{n}\dsum\limits_{j=0}^{y-1}(-1)^{j}B_{n}\left( x+%
\frac{j}{y}\right) ,  \label{MF2}
\end{equation}%
where $B_{n}(x)$ and $E_{n}(x)$ denote the Bernoulli polynomials and Euler
polynomials, respectively, cf. (\cite{Carlitz1952}, \cite%
{SrivastavaKurtSimsek}).
\end{remark}

To prove the multiplication formula of the generalized Apostol Bernoulli
polynomials, we need the following generating function which is defined by
Srivastava et al. \cite[pp. 254, Eq. (20)]{SrivastawaGargeSC}:

\begin{definition}
\label{DefBER}Let $a,b,c\in \mathbb{R}^{+}$ with $a\neq b,$ $x\in \mathbb{R}$
and $n\in \mathbb{N}_{0}$. Then the generalized Bernoulli polynomials $%
\mathfrak{B}_{n}^{(\alpha )}(x;\lambda ;a,b,c)$ of order $\alpha \in \mathbb{%
C}$ are defined by means of the following generating functions:%
\begin{equation}
f_{B}(x,a,b,c;\lambda ;\alpha )=\left( \frac{t}{\lambda b^{t}-a^{t}}\right)
^{\alpha }c^{xt}=\sum_{n=0}^{\infty }\mathfrak{B}_{n}^{(\alpha )}(x;\lambda
;a,b,c)\frac{t^{n}}{n!},  \label{1S}
\end{equation}%
where%
\begin{equation*}
\left\vert t\ln (\frac{a}{b})+\ln \lambda \right\vert <2\pi
\end{equation*}%
and%
\begin{equation*}
1^{\alpha }=1.
\end{equation*}
\end{definition}

Observe that if we set $\lambda =1$ in (\ref{1S}), we have%
\begin{equation}
\left( \frac{t}{b^{t}-a^{t}}\right) ^{\alpha }c^{xt}=\sum_{n=0}^{\infty }%
\mathfrak{B}_{n}^{(\alpha )}(x;a,b,c)\frac{t^{n}}{n!}.  \label{9}
\end{equation}%
If we set $x=0$ in (\ref{9}), we obtain%
\begin{equation}
\left( \frac{t}{b^{t}-a^{t}}\right) ^{\alpha }=\sum_{n=0}^{\infty }\mathfrak{%
B}_{n}^{(\alpha )}(a,b)\frac{t^{n}}{n!},  \label{8}
\end{equation}%
with of course, $\mathfrak{B}_{n}^{(\alpha )}(x;a,b,c)=\mathfrak{B}%
_{n}^{(\alpha )}(a,b)$, cf. (\cite{luo14}-\cite{lou15}, \cite{KimSimskJKM}, 
\cite{YsimsekKim}, \cite{kurtSimsek}, \cite{Mali}, \cite{OzdenAML}, \cite%
{OzdenSrivastava}, \cite{srivas18}, \cite{srivastava11}, \cite%
{Srivastava2011}, \cite{SrivastawaGargeSC}). If we set $\alpha =1$ in (\ref%
{8}) and (\ref{9}), we have%
\begin{equation}
\frac{t}{b^{t}-a^{t}}=\sum_{n=0}^{\infty }\mathfrak{B}_{n}(a,b)\frac{t^{n}}{%
n!}  \label{4}
\end{equation}%
and 
\begin{equation}
\left( \frac{t}{b^{t}-a^{t}}\right) c^{xt}=\sum_{n=0}^{\infty }\mathfrak{B}%
_{n}(x;a,b,c)\frac{t^{n}}{n!},  \label{5}
\end{equation}%
which have been studied by Luo et al. \cite{luo14}-\cite{lou15}. Moreover,
by substituting $a=1$ and $b=c=e$ into (\ref{1S}), then we arrive at the
Apostol-Bernoulli polynomials $\mathcal{B}_{n}(x;\lambda )$, which are
defined by means of the following generating function%
\begin{equation*}
\left( \frac{t}{\lambda e^{t}-1}\right) e^{xt}=\sum_{n=0}^{\infty }\mathcal{B%
}_{n}(x;\lambda )\frac{t^{n}}{n!},
\end{equation*}%
These polynomials $\mathcal{B}_{n}(x;\lambda )$ have been introduced and
investigated by many Mathematicians cf. (\cite{apostol}, \cite{KimChiARXIV}, 
\cite{KimkimJang}, \cite{KimSimskJKM}, \cite{burakSimsek}, \cite{luo13}, 
\cite{OzdenSrivastava}, \cite{SimsekJMAA}, \cite{srivas18}). When $a=\lambda
=1$ and $b=c=e$ into (\ref{4}) and (\ref{5}), $\mathfrak{B}_{n}(a,b)$ and $%
\mathfrak{B}_{n}(x;a,b,c)$ reduce to the classical Bernoulli numbers and the
classical Bernoulli polynomials, respectively, cf. \cite{AgohDilcher}-\cite%
{walum}.

\begin{remark}
The constraints on $\left\vert t\right\vert $, which we have used in
Definition \ref{DefBER} and (\ref{1Ssf}), are meant to ensure that the
generating function in (\ref{9})and (\ref{1Ssf}) are analytic throughout the
prescribed open disks in complex $t$-plane (centred at the origin $t=0$) in
order to have the corresponding convergent Taylor-Maclaurin series expansion
(about the origin $t=0$) occurring on the their right-hand side (with a
positive radius of convergence) cf. \cite{srivastava11}.
\end{remark}

\begin{theorem}
\label{T11a}Let $y\in \mathbb{N}$. Then we have%
\begin{equation*}
\mathfrak{B}_{n}(yx;\lambda
;a,b,b)=\dsum\limits_{l=0}^{n}\dsum\limits_{j=0}^{y-1}\left( 
\begin{array}{c}
n \\ 
l%
\end{array}%
\right) \lambda ^{j}y^{l-1}\left( (y-1-j)\ln a\right) ^{n-l}\mathfrak{B}%
_{l}\left( x+\frac{j}{y};\lambda ^{y};a,b,b\right) .
\end{equation*}
\end{theorem}

\begin{proof}
Substituting $c=b$ and $\alpha =1$ into (\ref{1S}), we get%
\begin{equation*}
\sum_{n=0}^{\infty }\mathfrak{B}_{n}(x;\lambda ;a,b,c)\frac{t^{n}}{n!}=\frac{%
1}{y}\sum_{j=0}^{y-1}\lambda ^{j}\frac{yt}{\lambda ^{y}b^{yt}-a^{yt}}%
b^{\left( \frac{x+j}{y}\right) yt}a^{t(y-j-1)}.
\end{equation*}%
Therefore%
\begin{eqnarray*}
&&\sum_{n=0}^{\infty }\mathfrak{B}_{n}(x;\lambda ;a,b,c)\frac{t^{n}}{n!} \\
&=&\sum_{n=0}^{\infty }\dsum\limits_{l=0}^{n}\dsum\limits_{j=0}^{y-1}\left( 
\begin{array}{c}
n \\ 
l%
\end{array}%
\right) \lambda ^{j}\left( (y-1-j)\ln a\right) ^{n-l}y^{l-1}\mathfrak{B}%
_{l}\left( \frac{x+j}{y};\lambda ^{y};a,b,b\right) \frac{t^{n}}{n!}.
\end{eqnarray*}%
Comparing the coefficients of $\frac{t^{n}}{n!}$ on both sides of the above
equation, we get%
\begin{equation*}
\mathfrak{B}_{n}^{(\alpha )}(x;\lambda
;a,b,c)=\dsum\limits_{l=0}^{n}\dsum\limits_{j=0}^{y-1}\left( 
\begin{array}{c}
n \\ 
k%
\end{array}%
\right) \lambda ^{j}\left( (k-1-j)\ln a\right) ^{n-l}y^{l-1}\mathfrak{B}%
_{l}\left( \frac{x+j}{k};\lambda ^{y};a,b,b\right) .
\end{equation*}%
By replacing $x$ by $yx$ on both sides of the above equation, we arrive at
the desired result.
\end{proof}

\begin{remark}
Kurt and Simsek \cite{kurtSimsek} proved multiplication formula for the
generalized Bernoulli polynomials of order $\alpha $. When $a=\lambda =1$
and $b=c=e$ into Theorem \ref{T11a}, we have the multiplication formula for
the Bernoulli polynomials given by%
\begin{equation}
B_{n}(yx)=y^{n-1}\dsum\limits_{j=0}^{y-1}B_{n}\left( x+\frac{j}{y}\right) ,
\label{MF1}
\end{equation}%
cf. (\cite{apostol}, \cite{Carlitz1952}, \cite{Carlitz}, \cite{DERE}, \cite%
{KimSimskJKM}, \cite{lou15}, \cite{luo13}, \cite{luo2003}, \cite{Mali}, \cite%
{LuoSrivatava2010}, \cite{srivas18}, \cite{SrivastawaGargeSC}).
\end{remark}

If $f$ is a \textit{normalized} polynomial such that it satisfies the formula%
\begin{equation}
f_{n}(yx)=y^{n-1}\dsum\limits_{j=0}^{y-1}f_{n}\left( x+\frac{j}{y}\right) ,
\label{w}
\end{equation}%
then $f$ is the $y$th degree Bernoulli polynomial due to (\ref{MF1}) cf. (%
\cite{Carlitz1952}, \cite{walum}). According to Nielsen \cite{Carlitz1952},
if a normalized polynomial satisfies (\ref{MF1}) for a single value of $y>1$%
, then it is identical with $B_{m}(x)$. Consequently, if a normalized
polynomial satisfies (\ref{mmF}) for a single value of $y>1$, then it is
identical with $\mathcal{H}_{n}(x;u;1,b,b;\lambda )$. The formula (\ref{MF2}%
) is different. Therefore, for $y$ is an even positive integer, Carlitz \cite%
[Eq-(1.4)]{Carlitz1952} considered the following equation:%
\begin{equation*}
g_{n-1}(yx)=-\frac{2y^{n-1}}{n}\dsum\limits_{j=0}^{y-1}(-1)^{j}f_{n}\left( x+%
\frac{j}{y}\right) ,
\end{equation*}%
where $g_{n-1}(x)$\ and $f_{n}(x)$ denote the normalized polynomials of
degree $n-1$ and $n$, respectively. More precisely, as Carlitz has pointed
out \cite[p. 184]{Carlitz1952}, if $y$ is a fixed even integer $\geq 2$ and $%
f_{n}(x)$\ is an arbitrary normalized polynomial of degree $n$, then (\ref%
{MF2}) determines $g_{n-1}(x)$\ as a normalized polynomial of degree $n-1$.
Thus, for a single value$\ y$, (\ref{MF2}) does not suffice to determine the
normalized polynomials $g_{n-1}(x)$\ and $f_{n}(x)$.

\begin{remark}
According to (\ref{w}), the set of normalized polynomials $\left\{
f_{n}(x)\right\} $ is an Appell set, cf. \cite{Carlitz1952}.
\end{remark}

We now modify (\ref{4ge1}) as follows:%
\begin{equation}
\frac{\left( a^{t}-\xi \right) c^{xt}}{\lambda b^{t}-\xi }%
=\sum_{n=0}^{\infty }\mathcal{H}_{n}(x;\xi ;a,b,c;\lambda )\frac{t^{n}}{n!}
\label{M1b}
\end{equation}%
where%
\begin{equation*}
\xi ^{r}=1,\text{ }\xi \neq 1.
\end{equation*}

The polynomial $\mathcal{H}_{n}(x;\xi ;a,b,c;\lambda )$ is a normalized
polynomial of degree $m$ in $x$. The polynomial \QTR{cal}{H}$_{n}(x;\xi
;1,e,e;1)$ may be called Eulerian polynomials with parameter $\xi $. In
particular we note that%
\begin{equation*}
\mathcal{H}_{n}(x;-1;1,e,e;1)=E_{n}(x)
\end{equation*}%
since for $a=\lambda =1$, $b=c=e$, equation (\ref{M1b}) reduces to the
generating function for the Euler polynomials.

By means of equation (\ref{mMF}), it is easy to verify the following
multiplication formulas:

If $y$ is an odd positive integer, then we have%
\begin{eqnarray}
\mathcal{H}_{n-1}(yx;\xi ;a,b,b;\lambda ) &=&\frac{y^{n-1}}{n}%
\dsum\limits_{j=0}^{y-1}\left( \frac{\lambda }{\xi }\right) ^{j}\mathfrak{B}%
_{n}\left( x+\frac{j}{y};b;\lambda ^{y}\right)   \label{M1c} \\
&&-\frac{1}{\xi n}\dsum\limits_{k=0}^{n}\dsum\limits_{j=0}^{y-1}\left( \frac{%
\lambda }{\xi }\right) ^{j}y^{k-1}(\ln a)^{n-k}\mathfrak{B}_{k}\left( x+%
\frac{j}{y};b;\lambda ^{y}\right) ,  \notag
\end{eqnarray}%
where%
\begin{equation*}
\mathcal{H}_{k}\left( x+\frac{j}{y};\xi ^{y};1,b,b;\lambda ^{y}\right) =%
\mathfrak{B}_{n}\left( x+\frac{j}{y};b;\lambda ^{y}\right) \text{.}
\end{equation*}%
If $y$ is an even positive integer, then we have%
\begin{eqnarray}
\mathcal{H}_{n}(yx;\xi ;a,b,b;\lambda ) &=&\frac{y^{n}}{2}%
\dsum\limits_{j=0}^{y-1}\left( \frac{\lambda }{\xi }\right) ^{j}\mathfrak{E}%
_{n}\left( x+\frac{j}{y};b;\lambda ^{y}\right)   \label{M1d} \\
&&-\frac{1}{2\xi }\dsum\limits_{k=0}^{n}\dsum\limits_{j=0}^{y-1}\left( \frac{%
\lambda }{\xi }\right) ^{j}y^{k}(\ln a)^{n-k}\mathfrak{E}_{k}\left( x+\frac{j%
}{y};b;\lambda ^{y}\right) ,  \notag
\end{eqnarray}%
where%
\begin{equation*}
\mathcal{H}_{k}\left( x+\frac{j}{y};\xi ^{y};1,b,b;\lambda ^{y}\right) =%
\mathfrak{E}_{n}\left( x+\frac{j}{y};b;\lambda ^{y}\right) \text{,}
\end{equation*}%
where $\mathfrak{E}_{n}(x;a,b,c)$ denotes the generalized Euler polynomials,
which are defined by means of the following generating function:%
\begin{equation*}
\left( \frac{t}{b^{t}-a^{t}}\right) c^{xt}=\sum_{n=0}^{\infty }\mathfrak{E}%
_{n}(x;a,b,c)\frac{t^{n}}{n!}
\end{equation*}%
cf. (\cite{luo14}-\cite{lou15}, \cite{Kim Jang}, \cite{kurtSimsek}, \cite%
{OzdenSrivastava}, \cite{srivas18}, \cite{srivastava11}, \cite%
{Srivastava2011}, \cite{SrivastawaGargeSC}).

\begin{remark}
If we set $a=\lambda =1$ and $b=e$, then (\ref{M1c}) and (\ref{M1d}) reduce
to the following multiplication formulas, respectively:%
\begin{equation*}
H_{n-1}(yx;\xi )=\frac{y^{n-1}}{n}\left( 1-\frac{1}{\xi }\right)
\dsum\limits_{j=0}^{y-1}\frac{1}{\xi ^{j}}B_{n}\left( x+\frac{j}{y}\right) 
\end{equation*}%
cf. \cite[Eq. (3.3)]{Carlitz1952} and%
\begin{equation*}
H_{n}(yx;\xi )=\frac{y^{n}}{2}\left( 1-\frac{1}{\xi }\right)
\dsum\limits_{j=0}^{y-1}\frac{1}{\xi ^{j}}E_{n}\left( x+\frac{j}{y}\right) .
\end{equation*}%
Let $f_{n}(x)$ and $g_{n}(x)$ be normalized polynomials in the usual way.
Carlitz \cite[Eq. (3.4)]{Carlitz1952} defined the following equation:%
\begin{equation*}
g_{n-1}(yx)=\frac{(1-\rho )y^{n-1}}{n}\dsum\limits_{j=0}^{y-1}\rho
^{j}f_{n}\left( x+\frac{j}{y}\right) ,
\end{equation*}%
where $\rho $ is a fixed primitive $r$th root of unity, $r>1$, $y\equiv 0(%
\func{mod}r)$.
\end{remark}

\begin{remark}
If we set $a=\lambda =1$, $b=c=e$ and $\xi =-1$, then (\ref{M1c}) and (\ref%
{M1d}) reduce to (\ref{MF2}) and (\ref{MF}).
\end{remark}

\begin{remark}
Walum \cite{walum} defined multiplication formula for periodic functions as
follows:%
\begin{equation}
\vartheta (y)f(yx)=\dsum\limits_{j(y)}f\left( x+\frac{j}{y}\right) ,
\label{w1}
\end{equation}%
where $f$ is periodic with period $1$ and $j(y)$ under the summation sign
indicates that $j$ runs through a complete system of residues $\func{mod}y$.

Formulas (\ref{w}), (\ref{w1}) and other multiplication formulas related to
periodic functions and normalized polynomials occur in Franel's formula, in
the theory of the Dedekind sums and Hardy-Berndt sums, in the theory of the
zeta functions and $L$-functions and in the theory of periodic bounded
variation, cf. (\cite{Berndt}, \cite{berndt2}, \cite{walum}).
\end{remark}

\subsection{Generalized Eulerian type numbers and polynomials attached to
Dirichlet character}

In this section, we construct generating function, related to nonnegative
real parameters, for the generalized Eulerian type numbers and polynomials
attached to Dirichlet character. We also give some properties of these
polynomials and numbers.

\begin{definition}
Let $\chi $ be the Dirichlet character of conductor $f\in \mathbb{N}$. Let $%
x\in \mathbb{R}$, $a,b\in \mathbb{R}^{+},$ $(a\neq b),$ $\lambda \in \mathbb{%
C}$ and $u\in \mathbb{C\diagdown }\left\{ 1\right\} $. The\ generalized
Eulerian type polynomials $\mathcal{H}_{n,\chi }(x;u;a,b,c;\lambda )$ are
defined by means of the following generating function:%
\begin{equation}
\mathcal{F}_{\lambda ,\chi }(t,x;u,a,b,c)=\dsum\limits_{j=0}^{f-1}\frac{%
\left( a^{ft}-u^{f}\right) \chi (j)u^{f-j-1}c^{\left( \frac{x+j}{f}\right)
ft}}{\lambda ^{f}b^{ft}-u^{f}}=\sum_{n=0}^{\infty }\mathcal{H}_{n,\chi
}(x;u;a,b,c;\lambda )\frac{t^{n}}{n!}  \label{4ge2}
\end{equation}

with, of course%
\begin{equation*}
\mathcal{H}_{n,\chi }(0;u;a,b,c;\lambda )=\mathcal{H}_{n,\chi
}(u;a,b,c;\lambda ),
\end{equation*}%
where $\mathcal{H}_{n,\chi }(u;a,b,c;\lambda )$ denotes generalized Eulerian
type numbers.
\end{definition}

\begin{remark}
In the special case when $a=\lambda =1$ and $b=c=e$, the generalized
Eulerian type polynomials $\mathcal{H}_{n,\chi }(x;u;a,b,c;\lambda )$\ are
reduced to the Frobenius Euler polynomials which are defined by means of the
following generating function:%
\begin{equation*}
\dsum\limits_{j=0}^{f-1}\frac{\left( 1-u^{f}\right) \chi
(j)u^{f-j-1}e^{\left( \frac{x+j}{f}\right) ft}}{e^{ft}-u^{f}}%
=\sum_{n=0}^{\infty }H_{n,\chi }(x;u)\frac{t^{n}}{n!},
\end{equation*}%
cf. (\cite{Tsumura}, \cite{KimSimskJKM}, \cite{YsimsekKim}, \cite{SimsekBKMS}%
, \cite{SimsekJNT}, \cite{srivas18}). Substituting $u=-1$ into the above
equation, we have generating function of the generalized Euler polynomials
attached to Dirichlet character with odd conductor:%
\begin{equation*}
2\dsum\limits_{j=0}^{f-1}\frac{\chi (j)(-1)^{j}e^{\left( \frac{x+j}{f}%
\right) ft}}{e^{ft}+1}=\sum_{n=0}^{\infty }E_{n,\chi }(x)\frac{t^{n}}{n!},
\end{equation*}%
cf. (\cite{Tsumura}, \cite{SimsekBKMS}, \cite{SimsekJNT}, \cite{srivas18}).
\end{remark}

Combining (\ref{4ge1}) and (\ref{4ge2}), we obtain the following functional
equation:%
\begin{equation*}
\mathcal{F}_{\lambda ,\chi }(t,x;u,a,b,c;)=\dsum\limits_{j=0}^{f-1}\chi
(j)u^{f-j-1}F_{\lambda ^{f}}(ft,\frac{x+j}{f};u^{f},a,b,c).
\end{equation*}

By using the above functional equation we arrive at the following Theorem:

\begin{theorem}
\begin{equation*}
\mathcal{H}_{n,\chi }(x;u;a,b,c;\lambda )=f^{n}\dsum\limits_{j=0}^{f-1}\chi
(j)u^{f-j-1}\mathcal{H}_{n}(\frac{x+j}{f};u^{f};a,b,c;\lambda ^{f}).
\end{equation*}
\end{theorem}

\begin{theorem}
\begin{equation*}
\mathcal{H}_{n,\chi }(x;u;a,b,c;\lambda )=\sum_{j=0}^{n}\left( 
\begin{array}{c}
n \\ 
j%
\end{array}%
\right) \left( x\ln c\right) ^{n-j}\mathcal{H}_{j,\chi }(u;a,b,c;\lambda ).
\end{equation*}
\end{theorem}

\begin{proof}
By using (\ref{4ge2}), we get%
\begin{equation*}
\sum_{n=0}^{\infty }\mathcal{H}_{n,\chi }(u;a,b,c;\lambda )\frac{t^{n}}{n!}%
\sum_{n=0}^{\infty }\left( x\ln c\right) ^{n}\frac{t^{n}}{n!}%
=\sum_{n=0}^{\infty }\mathcal{H}_{n,\chi }(x;u;a,b,c;\lambda )\frac{t^{n}}{n!%
}.
\end{equation*}%
From the above equation, we obtain%
\begin{equation*}
\sum_{n=0}^{\infty }\left( \sum_{j=0}^{n}\left( 
\begin{array}{c}
n \\ 
j%
\end{array}%
\right) \left( x\ln c\right) ^{n-j}\mathcal{H}_{j,\chi }(u;a,b,c;\lambda
\right) \frac{t^{n}}{n!}=\sum_{n=0}^{\infty }\mathcal{H}_{n,\chi
}(x;u;a,b,c;\lambda )\frac{t^{n}}{n!}.
\end{equation*}%
Comparing the coefficients of $\frac{t^{n}}{n!}$ on both sides of the above
equation, we arrive at the desired result.
\end{proof}

\subsection{\textbf{Recurrence relation for the }generalized Eulerian type
polynomials}

In this section we are going to differentiate (\ref{4ge1}) with respect to
the variable $t$ to derive a recurrence relation for the generalized
Eulerian type polynomials. Therefore, we obtain the following differential
equation:%
\begin{eqnarray*}
\frac{\partial }{\partial t}F_{\lambda }(t,x;u,a,b,c) &=&\left( \ln a\right)
F_{\lambda }(t,x;u,a,b,c)+\frac{\ln a}{t}f_{B}(x,1,b,c;\frac{\lambda }{u};1)
\\
&&-\frac{\ln \left( b^{\lambda }\right) }{ut}F_{\lambda
}(t,x;u,a,b,c)f_{B}(1,1,b,b;\frac{\lambda }{u};1) \\
&&+\ln \left( c^{x}\right) F_{\lambda }(t,x;u,a,b,c).
\end{eqnarray*}%
By using this equation, we obtain a recurrence relation for the generalized
Eulerian type polynomials by the following theorem:

\begin{theorem}
Let $n\in \mathbb{N}$. We have%
\begin{eqnarray*}
n\mathcal{H}_{n}(x;u;a,b,c;\lambda ) &=&\left( \ln a\right) \left( n\mathcal{%
H}_{n-1}(x;u;a,b,c;\lambda )+\mathfrak{B}_{n}(x;\frac{\lambda }{u}%
;a,b,c)\right) \\
&&-\frac{\lambda \ln b}{u}\sum_{j=0}^{n}\left( 
\begin{array}{c}
n \\ 
j%
\end{array}%
\right) \mathcal{H}_{j}(x;u;a,b,c;\lambda )\mathfrak{B}_{n-j}(1;\frac{%
\lambda }{u};1,b,b) \\
&&+\left( \ln \left( c^{nx}\right) \right) \mathcal{H}_{n-1}(x;u;a,b,c;%
\lambda ),
\end{eqnarray*}%
where $\mathfrak{B}_{n}(x;\lambda ;a,b,c)$ denotes the generalized Bernoulli
polynomials of order $1$.
\end{theorem}

\begin{remark}
When $a=\lambda =1$ and $b=c=e$, the recurrence relation for the generalized
Eulerian type polynomials is reduced to%
\begin{equation*}
nH_{n}(x;u)=nxH_{n-1}(x;u)-\frac{1}{u}\sum_{j=0}^{n}\left( 
\begin{array}{c}
n \\ 
j%
\end{array}%
\right) H_{j}(x;u)\mathcal{B}_{n-j}(1;\frac{1}{u}).
\end{equation*}
\end{remark}

\section{New identities involving families of polynomials}

In this section, we derive some new identities related to the generalized
Bernoulli polynomials and numbers of order $1$, the Eulerian type
polynomials and the generalized array type polynomials.

\begin{theorem}
\label{T11b} The following relationship holds true:%
\begin{equation*}
\mathfrak{B}_{n}(x;\lambda ;a,b,b)=\sum_{j=0}^{n}\left( 
\begin{array}{c}
n \\ 
j%
\end{array}%
\right) \mathcal{H}_{j}(x;\lambda ^{-1};a,\frac{b}{a},\frac{b}{a};1)%
\mathfrak{B}_{n-j}(x-1;\lambda ;1,a,a).
\end{equation*}
\end{theorem}

\begin{proof}
\begin{equation*}
\sum_{n=0}^{\infty }\mathfrak{B}_{n}(x;\lambda ;a,b,b)\frac{t^{n}}{n!}%
=\left( \frac{ta^{(x-1)t}}{\lambda a^{t}-1}\right) \left( \frac{\left(
a^{t}-\lambda ^{-1}\right) \left( \frac{b}{a}\right) ^{xt}}{\left( \frac{b}{a%
}\right) ^{t}-\lambda ^{-1}}\right) .
\end{equation*}%
Combining (\ref{1S}) and (\ref{4ge1}) with the above equation, we get%
\begin{equation*}
\sum_{n=0}^{\infty }\mathfrak{B}_{n}(x;\lambda ;a,b,b)\frac{t^{n}}{n!}%
=\sum_{n=0}^{\infty }\mathfrak{B}_{n}(x-1;\lambda ;1,a,a)\frac{t^{n}}{n!}%
\sum_{n=0}^{\infty }\mathcal{H}_{n}(x;\lambda ^{-1};a,\frac{b}{a},\frac{b}{a}%
;1)\frac{t^{n}}{n!}.
\end{equation*}%
Therefore%
\begin{equation*}
\sum_{n=0}^{\infty }\mathfrak{B}_{n}(x;\lambda ;a,b,b)\frac{t^{n}}{n!}%
=\sum_{n=0}^{\infty }\left( \sum_{j=0}^{n}\left( 
\begin{array}{c}
n \\ 
j%
\end{array}%
\right) \mathcal{H}_{j}(x;\lambda ^{-1};a,\frac{b}{a},\frac{b}{a};1)%
\mathfrak{B}_{n-j}(x-1;\lambda ;1,a,a)\right) \frac{t^{n}}{n!}.
\end{equation*}%
Comparing the coefficients of $\frac{t^{n}}{n!}$ on both sides of the above
equation, we arrive at the desired result.
\end{proof}

Relationship between the generalized Bernoulli numbers and the Frobenius
Euler numbers is given by the following result:

\begin{theorem}
The following relationship holds true:%
\begin{equation}
\mathfrak{B}_{n}(\lambda ;a,b)=\frac{1}{\lambda -1}\sum_{j=0}^{n}\left( 
\begin{array}{c}
n \\ 
j%
\end{array}%
\right) j\left( \ln a^{-1}\right) ^{n-j}\left( \ln \left( \frac{b}{a}\right)
\right) ^{j}H_{j-1}\left( \lambda ^{-1}\right) .  \label{1Ssd}
\end{equation}
\end{theorem}

\begin{proof}
By using (\ref{1S}), we obtain%
\begin{equation*}
\sum_{n=0}^{\infty }\mathfrak{B}_{n}(\lambda ;a,b)\frac{t^{n}}{n!}=\frac{%
ta^{-t}}{\lambda -1}\left( \frac{1-\lambda ^{-1}}{e^{t\ln \left( \frac{b}{a}%
\right) }-\lambda ^{-1}}\right) .
\end{equation*}

From the above equation, we get%
\begin{equation*}
\sum_{n=0}^{\infty }\mathfrak{B}_{n}(\lambda ;a,b)\frac{t^{n}}{n!}=\frac{1}{%
\lambda -1}\sum_{n=0}^{\infty }\left( \ln \left( \frac{1}{a}\right) \right)
^{n}\frac{t^{n}}{n!}\sum_{n=0}^{\infty }n\mathcal{H}_{n}(\lambda
^{-1})\left( \ln \left( \frac{b}{a}\right) \right) ^{n}\frac{t^{n}}{n!}.
\end{equation*}%
Therefore%
\begin{equation*}
\sum_{n=0}^{\infty }\mathfrak{B}_{n}(\lambda ;a,b)\frac{t^{n}}{n!}%
=\sum_{n=0}^{\infty }\left( \sum_{j=0}^{n}\left( 
\begin{array}{c}
n \\ 
j%
\end{array}%
\right) \frac{j\left( \ln a^{-1}\right) ^{n-j}\left( \ln \left( \frac{b}{a}%
\right) \right) ^{j}}{\lambda -1}H_{j-1}\left( \lambda ^{-1}\right) \right) 
\frac{t^{n}}{n!}.
\end{equation*}

Comparing the coefficients of $\frac{t^{n}}{n!}$ on both sides of the above
equation, we arrive at the desired result.
\end{proof}

\begin{remark}
By substituting $a=1$ and $b=e$ into (\ref{1Ssd}), we have%
\begin{equation*}
\mathcal{B}_{n}(\lambda )=\frac{n}{\lambda -1}H_{n-1}(\lambda ^{-1}),
\end{equation*}%
cf. \cite{KimSimskJKM}.
\end{remark}

Relationship between the generalized Eulerian type polynomials and
generalized array type polynomials are given by the following theorem:

\begin{theorem}
The following relationship holds true:%
\begin{equation*}
\mathcal{H}_{n}(x;u;a,b,b;\lambda )=\sum_{k=0}^{\infty }\sum_{m=0}^{\infty
}\sum_{d=0}^{n}\left( 
\begin{array}{c}
m+k-1 \\ 
m%
\end{array}%
\right) \left( 
\begin{array}{c}
n \\ 
d%
\end{array}%
\right) \frac{k!\left( \ln a^{m}\right) ^{n-d}}{u^{m+k}}\mathcal{S}%
_{k}^{d}(x;a,b;\lambda ).
\end{equation*}
\end{theorem}

\begin{proof}
From (\ref{4ge1}), we obtain%
\begin{equation*}
\sum_{n=0}^{\infty }\mathcal{H}_{n}(x;u;a,b,c;\lambda )\frac{t^{n}}{n!}%
=\sum_{k=0}^{\infty }\left( \frac{\lambda b^{t}-a^{t}}{u-a^{t}}\right)
^{k}b^{xt}.
\end{equation*}%
Combining (\ref{ab0}) with the above equation, we get%
\begin{equation*}
\sum_{n=0}^{\infty }\mathcal{H}_{n}(x;u;a,b,b;\lambda )\frac{t^{n}}{n!}%
=\sum_{k=0}^{\infty }\frac{k!}{\left( u-a^{t}\right) ^{k}}\sum_{n=0}^{\infty
}\mathcal{S}_{k}^{n}(x;a,b;\lambda )\frac{t^{n}}{n!}.
\end{equation*}%
From the above equation, we get%
\begin{equation*}
\sum_{n=0}^{\infty }\mathcal{H}_{n}(x;u;a,b,b;\lambda )\frac{t^{n}}{n!}%
=\sum_{n=0}^{\infty }\sum_{k=0}^{\infty }\frac{k!\mathcal{S}%
_{k}^{n}(x;a,b;\lambda )}{u^{k}\left( 1-\frac{a^{t}}{u}\right) ^{k}}\frac{%
t^{n}}{n!}\text{.}
\end{equation*}%
Now we assume $\left\vert \frac{a^{t}}{u}\right\vert <1$ in the above
equation; thus we get%
\begin{eqnarray*}
&&\sum_{n=0}^{\infty }\mathcal{H}_{n}(x;u;a,b,b;\lambda )\frac{t^{n}}{n!} \\
&=&\sum_{n=0}^{\infty }\sum_{k=0}^{\infty }\sum_{m=0}^{\infty }\left( 
\begin{array}{c}
m+k-1 \\ 
m%
\end{array}%
\right) \frac{k!\mathcal{S}_{k}^{n}(x;a,b;\lambda )}{u^{k+m}}\frac{%
a^{mt}t^{n}}{n!}.
\end{eqnarray*}%
Therefore%
\begin{eqnarray*}
&&\sum_{n=0}^{\infty }\mathcal{H}_{n}(x;u;a,b,b;\lambda )\frac{t^{n}}{n!} \\
&=&\sum_{n=0}^{\infty }\left( \sum_{k=0}^{\infty }\sum_{m=0}^{\infty
}\sum_{d=0}^{n}\left( 
\begin{array}{c}
m+k-1 \\ 
m%
\end{array}%
\right) \left( 
\begin{array}{c}
n \\ 
d%
\end{array}%
\right) \frac{k!\left( \ln a^{m}\right) ^{n-d}}{u^{m+k}}\mathcal{S}%
_{k}^{d}(x;a,b;\lambda ).\right) \frac{t^{n}}{n!}.
\end{eqnarray*}%
Comparing the coefficients of $\frac{t^{n}}{n!}$ on both sides of the above
equation, we arrive at the desired result.
\end{proof}

\begin{remark}
Substituting $a=1$ into the above Theorem and noting that $d=n$, we deduce
the following identity:%
\begin{equation*}
\mathcal{H}_{n}(x;u;1,b,b;\lambda )=\sum_{k=0}^{\infty }\frac{k!}{(u-1)^{k}}%
\mathcal{S}_{k}^{n}(x;1,b;\lambda )
\end{equation*}%
which upon setting $\lambda =1$ and $b=e$, yields%
\begin{equation*}
H_{n}(x;u)=\sum_{k=0}^{n}\frac{k!}{(u-1)^{k}}\mathcal{S}_{k}^{n}(x)
\end{equation*}%
which was found by Chang and Ha \cite[Lemma 1]{Chan}.
\end{remark}

\section{Relationship between the generalized Bernoulli polynomials and the
generalized array type polynomials}

In this section, we give some applications related to the generalized
Bernoulli polynomials, generalized array type polynomials. We derive many
identities involving these polynomials. By using same method with Agoh and
Dilcher's \cite{AgohDilcher}, we give the following Theorem:

\begin{theorem}
\label{L-AD}%
\begin{equation}
\left( \frac{\lambda b^{t}-a^{t}}{t}\right)
^{k}b^{xt}=\dsum\limits_{n=0}^{\infty }\frac{\mathcal{S}_{k}^{n+k}(x;a,b;%
\lambda )}{\binom{n+k}{k}}\frac{t^{n}}{n!}.  \label{w3}
\end{equation}
\end{theorem}

\begin{proof}
Combining\ (\ref{ab0}) and (\ref{ab1}), we get 
\begin{eqnarray*}
\left( \frac{\lambda b^{t}-a^{t}}{t}\right) ^{k}b^{xt} &=&\frac{1}{t^{k}}%
\dsum\limits_{n=0}^{\infty }\frac{k!}{n!}S_{k}^{n}\left( x,a,b;\lambda
\right) t^{n} \\
&=&\dsum\limits_{n=0}^{\infty }\frac{k!}{n!}S_{k}^{n+k}\left( x,a,b;\lambda
\right) t^{n-k}.
\end{eqnarray*}%
From the above equation, we arrive at the desired result.
\end{proof}

\begin{remark}
By setting $x=0$, $a=\lambda =1$ and $b=e$, Theorem \ref{L-AD} yields the
corresponding result which is proven by Agoh and Dilcher \cite{AgohDilcher}.
\end{remark}

\begin{theorem}
\label{Tw3}%
\begin{eqnarray*}
&&\left( n+k\right) \frac{\mathcal{S}_{k}^{n+k}(x;a,b;\lambda )}{\binom{n+k}{%
k}}-xn\frac{\mathcal{S}_{k}^{n+k-1}(x;a,b;\lambda )}{\binom{n+k-1}{k}} \\
&=&\sum_{j=0}^{n}\frac{\left( 
\begin{array}{c}
n \\ 
j%
\end{array}%
\right) }{\left( 
\begin{array}{c}
j+k-1 \\ 
k-1%
\end{array}%
\right) }\mathcal{S}_{k-1}^{j+k-1}(x;a,b;\lambda )\left( \ln \left(
b^{\lambda k}\right) \left( \ln (b)\right) ^{n-j}-\ln \left( a^{k}\right)
\left( \ln (a)\right) ^{n-j}\right) .
\end{eqnarray*}
\end{theorem}

\begin{proof}
By differentiating both sides of equation (\ref{w3}) with respect to the
variable $t$, after some elementary calculations, we get the formula
asserted by Theorem \ref{Tw3}.
\end{proof}

\begin{theorem}
\label{Teo19} The following relationship holds true:%
\begin{equation*}
\mathcal{S}_{k-1}^{n}(x+y;a,b;\lambda )=\sum_{j=0}^{n}\frac{\left( 
\begin{array}{c}
n \\ 
j%
\end{array}%
\right) \left( 
\begin{array}{c}
n+k-1 \\ 
k-1%
\end{array}%
\right) }{\left( 
\begin{array}{c}
j+k \\ 
k%
\end{array}%
\right) }\mathcal{S}_{k}^{j+k}(x;a,b;\lambda )\mathfrak{B}_{n-j}(y;\lambda
;a,b,b).
\end{equation*}
\end{theorem}

\begin{proof}
We set%
\begin{equation*}
\left( \frac{\lambda b^{t}-a^{t}}{t}\right) ^{k}b^{xt}\left( \frac{tb^{yt}}{%
\lambda b^{t}-a^{t}}\right) =\left( \frac{\lambda b^{t}-a^{t}}{t}\right)
^{k-1}b^{(x+y)t}.
\end{equation*}%
Combining (\ref{w3}) and (\ref{5}) with the above equation, we get%
\begin{equation*}
\sum_{n=0}^{\infty }\frac{\mathcal{S}_{k-1}^{n+k-1}(x+y;a,b;\lambda )}{%
\binom{n+k-1}{k-1}}\frac{t^{n}}{n!}=\dsum\limits_{n=0}^{\infty }\mathfrak{B}%
_{n}(y;\lambda ;a,b,b)\frac{t^{n}}{n!}\sum_{n=0}^{\infty }\frac{\mathcal{S}%
_{k}^{n+k}(x;a,b;\lambda )}{\binom{n+k}{k}}\frac{t^{n}}{n!}.
\end{equation*}

Therefore%
\begin{equation*}
\sum_{n=0}^{\infty }\frac{\mathcal{S}_{k-1}^{n+k-1}(x+y;a,b;\lambda )}{%
\binom{n+k-1}{k-1}}\frac{t^{n}}{n!}=\sum_{n=0}^{\infty }\left( \sum_{j=0}^{n}%
\frac{\left( 
\begin{array}{c}
n \\ 
j%
\end{array}%
\right) }{\left( 
\begin{array}{c}
j+k \\ 
k%
\end{array}%
\right) }\mathcal{S}_{k}^{j+k}(x;a,b;\lambda )\mathfrak{B}_{n-j}(y;\lambda
;a,b,b)\right) \frac{t^{n}}{n!}.
\end{equation*}%
Comparing the coefficients of $\frac{t^{n}}{n!}$ on both sides of the above
equation, we arrive at the desired result.
\end{proof}

\begin{remark}
By setting $x=y=0$, $a=\lambda =1$ and $b=e$, Theorem \ref{Teo19} yields the
corresponding result which is proven by Agoh and Dilcher \cite{AgohDilcher}.
\end{remark}

\begin{theorem}
The following relationship holds true:%
\begin{equation*}
\mathfrak{B}_{n}^{(u-v)}(x+y;\lambda ;a,b,b)=\sum_{j=0}^{n}\frac{\left( 
\begin{array}{c}
n \\ 
j%
\end{array}%
\right) }{\left( 
\begin{array}{c}
n+v \\ 
v%
\end{array}%
\right) }\mathcal{S}_{v}^{j+v}(x;a,b;\lambda )\mathfrak{B}%
_{n-j}^{(u)}(y;\lambda ;a,b,b).
\end{equation*}
\end{theorem}

\begin{proof}
We set%
\begin{equation}
\left( \frac{\lambda b^{t}-a^{t}}{t}\right) ^{v}b^{xt}\left( \frac{t}{%
\lambda b^{t}-a^{t}}\right) ^{u}b^{yt}=\left( \frac{t}{\lambda b^{t}-a^{t}}%
\right) ^{u-v}b^{\left( x+y\right) t}.  \label{wc3}
\end{equation}%
Combining (\ref{w3}) and (\ref{1S}) with the above equation, by using same
calculations with the proof of Theorem \ref{Teo19}, we arrive at the desired
result.
\end{proof}

\section{Application of the Laplace transform to the generating functions
for the generalized Bernoulli polynomials and the generalized array type
polynomials}

In this section, we give an application of the Laplace transform to the
generating function for the generalized Bernoulli polynomials and the
generalized array type polynomials. We obtain interesting series
representation for the families of these polynomials.

By using (\ref{wc3}), we obtain%
\begin{eqnarray*}
&&\dsum\limits_{n=0}^{\infty }\mathfrak{B}_{n}^{(u-v)}(\lambda ;a,b,b)\frac{%
t^{n}}{n!}e^{-t(y-x)\ln b} \\
&=&\dsum\limits_{n=0}^{\infty }\left( \sum_{j=0}^{n}\frac{\left( 
\begin{array}{c}
n \\ 
j%
\end{array}%
\right) }{\left( 
\begin{array}{c}
n+v \\ 
v%
\end{array}%
\right) }\mathcal{S}_{v}^{j+v}(x;a,b;\lambda )\mathfrak{B}%
_{n-j}^{(u)}(\lambda ;a,b,b)\right) \frac{t^{n}}{n!}e^{-ty\ln b}.
\end{eqnarray*}%
Integrate this equation (by parts) with respect to $t$ from $0$ to $\infty $%
, we get%
\begin{eqnarray*}
&&\dsum\limits_{n=0}^{\infty }\frac{\mathfrak{B}_{n}^{(u-v)}(\lambda ;a,b,b)%
}{n!}\dint\limits_{0}^{\infty }t^{n}e^{-t(y-x)\ln b}dt \\
&=&\dsum\limits_{n=0}^{\infty }\left( \frac{1}{n!}\sum_{j=0}^{n}\frac{\left( 
\begin{array}{c}
n \\ 
j%
\end{array}%
\right) }{\left( 
\begin{array}{c}
n+v \\ 
v%
\end{array}%
\right) }\mathcal{S}_{v}^{j+v}(x;a,b;\lambda )\mathfrak{B}%
_{n-j}^{(u)}(\lambda ;a,b,b)\right) \dint\limits_{0}^{\infty }t^{n}e^{-ty\ln
b}dt.
\end{eqnarray*}%
By using Laplace transform in the above equation, we arrive at the following
Theorem:

\begin{theorem}
\label{TeoL}The following relationship holds true:%
\begin{equation*}
\dsum\limits_{n=0}^{\infty }\frac{\mathfrak{B}_{n}^{(u-v)}(\lambda ;a,b,b)}{%
(\ln b^{y-x})^{n+1}}=\dsum\limits_{n=0}^{\infty }\sum_{j=0}^{n}\frac{\left( 
\begin{array}{c}
n \\ 
j%
\end{array}%
\right) }{\left( 
\begin{array}{c}
n+v \\ 
v%
\end{array}%
\right) }\frac{\mathcal{S}_{v}^{j+v}(x;a,b;\lambda )\mathfrak{B}%
_{n-j}^{(u)}(\lambda ;a,b,b)}{\left( \ln b^{y}\right) ^{n+1}}.
\end{equation*}
\end{theorem}

\begin{remark}
When $a=\lambda =1$ and $b=e$, Theorem \ref{TeoL}\ is reduced to the
following result:%
\begin{equation*}
\dsum\limits_{n=0}^{\infty }\frac{B_{n}^{(u-v)}}{(y-x)^{n+1}}%
=\dsum\limits_{n=0}^{\infty }\sum_{j=0}^{n}\frac{\left( 
\begin{array}{c}
n \\ 
j%
\end{array}%
\right) }{\left( 
\begin{array}{c}
n+v \\ 
v%
\end{array}%
\right) }\frac{S_{v}^{j+v}(x)B_{n-j}^{(u)}}{y^{n+1}}.
\end{equation*}
\end{remark}

\section{Applications the $p$-adic integral to the family of the normalized
polynomials and the generalized $\protect\lambda $-Stirling type numbers}

By using the $p$-adic integrals on $%
\mathbb{Z}
_{p}$, we derive some new identities related to the Bernoulli numbers, the
Euler numbers, the generalized Eulerian type numbers and the generalized $%
\lambda $-Stirling type numbers.

In order to prove the main results in this section, we recall each of the
following known results related to the $p$-adic integral.

Let $p$ be a fixed prime. It is known that%
\begin{equation*}
\mu _{q}(x+p^{N}\mathbb{Z}_{p})=\frac{q^{x}}{\left[ p^{N}\right] _{q}}
\end{equation*}%
is a distribution on $\mathbb{Z}_{p}$ for $q\in \mathbb{C}_{p}$ with $\mid
1-q\mid _{p}<1$, cf. \cite{T. Kim}. Let $UD\left( \mathbb{Z}_{p}\right) $ be
the set of uniformly differentiable functions on $\mathbb{Z}_{p}$. The $p$%
-adic $q$-integral of the function $f\in UD\left( \mathbb{Z}_{p}\right) $ is
defined by Kim \cite{T. Kim} as follows:%
\begin{equation*}
\int_{\mathbb{Z}_{p}}f(x)d\mu _{q}(x)=\lim_{N\rightarrow \infty }\frac{1}{%
[p^{N}]_{q}}\sum_{x=0}^{p^{N}-1}f(x)q^{x},
\end{equation*}%
where%
\begin{equation*}
\left[ x\right] =\frac{1-q^{x}}{1-q}.
\end{equation*}%
From this equation, the\textit{\ bosonic} $p$-adic integral ($p$-adic
Volkenborn integral) was considered from a physical point of view to the
bosonic limit $q\rightarrow 1$, as follows (\cite{T. Kim}):%
\begin{equation}
\int\limits_{\mathbb{Z}_{p}}f\left( x\right) d\mu _{1}\left( x\right) =%
\underset{N\rightarrow \infty }{\lim }\frac{1}{p^{N}}\sum_{x=0}^{p^{N}-1}f%
\left( x\right) ,  \label{M}
\end{equation}%
where%
\begin{equation*}
\mu _{1}\left( x+p^{N}\mathbb{Z}_{p}\right) =\frac{1}{p^{N}}.
\end{equation*}%
The $p$-adic $q$-integral is used in many branch of mathematics,
mathematical physics and other areas cf. (\cite{Amice}, \cite{T. Kim}, \cite%
{KimSimskJKM}, \cite{Schikof}, \cite{K. Shiratani}, \cite{SimsekJMAA}, \cite%
{SimsekADEA}, \cite{srivas18}, \cite{Volkenborn}).

By using (\ref{M}), we have the Witt's formula for the Bernoulli numbers $%
B_{n}$ as follows:%
\begin{equation}
\int\limits_{\mathbb{Z}_{p}}x^{n}d\mu _{1}\left( x\right) =B_{n}  \label{M1}
\end{equation}%
cf. (\cite{Amice}, \cite{T. Kim}, \cite{Kim2006TMIC}, \cite%
{KimmskimlcjangJIA}, \cite{Schikof}, \cite{Volkenborn}).

We consider the \textit{fermionic} integral in contrast to the convential
bosonic, which is called the fermionic $p$-adic integral on $\mathbb{Z}_{p}$
cf. \cite{Kim2006TMIC}. That is%
\begin{equation}
\int\limits_{\mathbb{Z}_{p}}f\left( x\right) d\mu _{-1}\left( x\right) =%
\underset{N\rightarrow \infty }{\lim }\sum_{x=0}^{p^{N}-1}\left( -1\right)
^{x}f\left( x\right)   \label{Mm}
\end{equation}%
where%
\begin{equation*}
\mu _{1}\left( x+p^{N}\mathbb{Z}_{p}\right) =\frac{(-1)^{x}}{p^{N}}
\end{equation*}%
cf. \cite{Kim2006TMIC}. By using (\ref{Mm}), we have the Witt's formula for
the Euler numbers $E_{n}$ as follows:%
\begin{equation}
\int\limits_{\mathbb{Z}_{p}}x^{n}d\mu _{-1}\left( x\right) =E_{n},
\label{Mm1}
\end{equation}%
cf. (\cite{Kim2006TMIC}, \cite{KimmskimlcjangJIA}, \cite{SimsekADEA}, \cite%
{srivas18}).

The Volkenborn integral in terms of the Mahler coefficients is given by the
following Theorem:

\begin{theorem}
\label{TSHiR}Let%
\begin{equation*}
f(x)=\sum_{j=0}^{\infty }a_{j}\left( 
\begin{array}{c}
x \\ 
j%
\end{array}%
\right) \in UD\left( \mathbb{Z}_{p}\right) .
\end{equation*}%
Then%
\begin{equation*}
\int\limits_{\mathbb{Z}_{p}}f(x)d\mu _{1}\left( x\right) =\sum_{j=0}^{\infty
}a_{j}\frac{(-1)^{j}}{j+1}.
\end{equation*}
\end{theorem}

Proof of Theorem \ref{TSHiR} was given by Schikhof \cite{Schikof}.

\begin{theorem}
\label{L1}%
\begin{equation*}
\int\limits_{\mathbb{Z}_{p}}\left( 
\begin{array}{c}
x \\ 
j%
\end{array}%
\right) d\mu _{1}\left( x\right) =\frac{(-1)^{j}}{j+1}.
\end{equation*}%
Proof of Theorem \ref{L1} was given by Schikhof \cite{Schikof}.
\end{theorem}

\begin{theorem}
The following relationship holds true:%
\begin{equation}
B_{m}=\frac{1}{\ln ^{m}b}\sum_{j=0}^{m}(-1)^{j}\frac{j!}{j+1}\mathcal{S}%
(m,j;1,b;1).  \label{1Ss}
\end{equation}
\end{theorem}

\begin{proof}
If we substitute $a=\lambda =1$ in Theorem \ref{T3}, we have%
\begin{equation*}
\left( \ln b^{x}\right) ^{m}=\sum_{j=0}^{m}\left( 
\begin{array}{c}
x \\ 
j%
\end{array}%
\right) j!\mathcal{S}(m,j;1,b;1).
\end{equation*}%
By applying the $p$-adic Volkenborn integral with Theorem \ref{L1}\ to the
both sides of the above equation, we arrive at the desired result.
\end{proof}

\begin{remark}
By substituting $b=1$ into (\ref{1Ss}), we have 
\begin{equation*}
B_{m}=\sum_{j=0}^{m}(-1)^{j}\frac{j!}{j+1}S(m,j)
\end{equation*}%
where $S(m,j)$\ denotes the Stirling numbers of the second kind cf. (\cite%
{ChanManna}, \cite{http}, \cite{KimChiARXIV}).
\end{remark}

\begin{theorem}
\label{Teo14} The following relationship holds true:%
\begin{eqnarray*}
&&\sum_{j=0}^{n}\binom{n}{j}\left( \ln a\right) ^{n-j}\left( \ln c\right)
^{j}B_{j}-u(\ln c)^{n}B_{n} \\
&=&\sum_{j=0}^{n}\binom{n}{j}\left( \ln c\right) ^{j}\left( \lambda \left( 
\mathcal{H}(u;a,b,c;\lambda )+\ln b\right) ^{n-j}-u\mathcal{H}%
_{n-j}(u;a,b,c;\lambda )\right) B_{j}.
\end{eqnarray*}
\end{theorem}

\begin{proof}
By using Theorem \ref{t8}, we have%
\begin{eqnarray}
&&\sum_{j=0}^{n}\binom{n}{j}\left( \ln a\right) ^{n-j}\left( \ln c\right)
^{j}x^{j}-u(\ln c)^{n}x^{n}  \label{M1a} \\
&=&\sum_{j=0}^{n}\binom{n}{j}\left( \ln c\right) ^{j}x^{j}\left( \lambda
\left( \mathcal{H}(u;a,b,c;\lambda )+\ln b\right) ^{n-j}-u\mathcal{H}%
_{n-j}(u;a,b,c;\lambda )\right) .  \notag
\end{eqnarray}%
By applying Volkenborn integral in (\ref{M}) to the both sides of the above
equation, we get%
\begin{eqnarray*}
&&\sum_{j=0}^{n}\binom{n}{j}\left( \ln a\right) ^{n-j}\left( \ln c\right)
^{j}\int\limits_{\mathbb{Z}_{p}}x^{j}d\mu (x)-u(\ln c)^{n}\int\limits_{%
\mathbb{Z}_{p}}x^{n}d\mu (x) \\
&=&\sum_{j=0}^{n}\binom{n}{j}\left( \ln c\right) ^{j}\left( \lambda \left( 
\mathcal{H}(u;a,b,c;\lambda )+\ln b\right) ^{n-j}-u\mathcal{H}%
_{n-j}(u;a,b,c;\lambda )\right) \int\limits_{\mathbb{Z}_{p}}x^{j}d\mu (x).
\end{eqnarray*}%
By substituting (\ref{M1}) into the above equation, we easily arrive at the
desired result.
\end{proof}

\begin{remark}
By substituting $b=c=e$ and $a=\lambda =1$ into Theorem \ref{Teo14}, we
arrive at the following nice identity:%
\begin{equation*}
B_{n}=\frac{1}{1-u}\sum_{j=0}^{n}\binom{n}{j}\left( \left( H(u)+1\right)
^{n-j}-uH_{n-j}(u)\right) B_{j}.
\end{equation*}
\end{remark}

\begin{theorem}
\label{Te15} The following relationship holds true:%
\begin{eqnarray*}
&&\sum_{j=0}^{n}\binom{n}{j}\left( \ln a\right) ^{n-j}\left( \ln c\right)
^{j}E_{j}-u(\ln c)^{n}E_{n} \\
&=&\sum_{j=0}^{n}\binom{n}{j}\left( \ln c\right) ^{j}\left( \lambda \left( 
\mathcal{H}(u;a,b,c;\lambda )+\ln b\right) ^{n-j}-u\mathcal{H}%
_{n-j}(u;a,b,c;\lambda )\right) E_{j}.
\end{eqnarray*}
\end{theorem}

\begin{proof}
Proof of Theorem \ref{Te15} is same as that of Theorem \ref{Teo14}.
Combining (\ref{Mm}), (\ref{M1a}) and (\ref{Mm1}), we easily arrive at the
desired result.
\end{proof}

\begin{remark}
By substituting $b=c=e$ and $a=\lambda =1$ into Theorem \ref{Te15}, we
arrive at the following nice identity:%
\begin{equation*}
E_{n}=\frac{1}{1-u}\sum_{j=0}^{n}\binom{n}{j}\left( \left( H(u)+1\right)
^{n-j}-uH_{n-j}(u)\right) E_{j}.
\end{equation*}
\end{remark}

\begin{acknowledgement}
The present investigation was supported by the \textit{Scientific Research
Project Administration of Akdeniz University.}
\end{acknowledgement}

\end{document}